\documentclass{svjour3}
\smartqed 
\usepackage{marginnote}
\usepackage{graphicx}
\usepackage{subfigure}
\usepackage{xspace}
\usepackage[T1]{fontenc}
\usepackage{comment}
\usepackage{graphicx}
\usepackage[headings]{fullpage}
\usepackage{amsmath, amssymb, amsfonts}
\usepackage{commath}
\usepackage{times}
\usepackage[usenames]{color}
\usepackage{enumerate}
\usepackage{theorem}
\mathchardef\mhyphen="2D 
\usepackage[colorlinks=true]{hyperref}

\usepackage{hyperref}
\usepackage{verbatim} 
\usepackage{algpseudocode}

\usepackage{ulem}
\theorembodyfont{\upshape}
\newcommand{\sbt}{\vcenter{\hbox{\tiny$\bullet$}}}
\newcommand{\Psidot}{\overset{\sbt}{\Psi}}
\newcommand{\threehalf}{{\textstyle \frac{3}{2}}}
\newcommand{\forth}{{\textstyle \frac{1}{4}}}

\newtheorem{algorithm}{Algorithm}

\newcommand {\half} {\mbox{$\frac{1}{2}$}}

\begin{document}

\title{Fast evaluation of far-field signals for time-domain wave propagation}

\author{Scott E. Field and Stephen R. Lau}

\institute{Scott E. Field \at
Department of Physics, Joint Space Science Institute,\\ 
and Maryland Center for Fundamental Physics\\
University of Maryland\\
College Park, MD 20742, USA\\[5pt]
Center for Radiophysics and Space Research,\\
Cornell University,\\
Ithaca, New York 14853, USA
\and
Stephen R. Lau \at
Department of Mathematics and Statistics\\
University of New Mexico\\
Albuquerque, NM 87131, USA
}
\maketitle 
\begin{abstract}
Time-domain simulation of wave phenomena on a finite computational 
domain often requires a fictitious outer boundary. An important 
practical issue is the specification of appropriate boundary 
conditions on this boundary, often conditions of complete transparency. 
Attention to this issue has been paid elsewhere, and here we consider 
a different, although related, issue: far-field signal recovery. 
Namely, from smooth data recorded on the outer boundary we 
wish to recover the far-field signal which would 
reach arbitrarily large distances. These signals encode 
information about interior scatterers and often correspond to actual 
measurements. This article expresses far-field signal recovery in 
terms of time-domain convolutions, each between a solution multipole
moment recorded at the boundary and a sum-of-exponentials kernel.
Each exponential corresponds to a pole term in the Laplace transform
of the kernel, a finite sum of simple poles. Greengard,
Hagstrom, and Jiang have derived the large-$\ell$ (spherical-harmonic
index) asymptotic expansion for the pole residues,
and their analysis shows that, when expressed in terms of the exact
sum-of-exponentials, large-$\ell$ 
signal recovery is plagued by cancellation errors.
Nevertheless, through an alternative
integral representation of the kernel and its subsequent
approximation by a {\em smaller}
number of exponential terms (kernel compression),
we are able to alleviate these errors and achieve accurate signal
recovery. We empirically examine scaling relations
between the parameters which determine a compressed kernel,
and perform numerical tests of signal "teleportation"
from one radial value $r_1$ to another $r_2$, including
the case $r_2=\infty$. We conclude with
a brief discussion on application to other hyperbolic equations
posed on non-flat geometries where waves undergo backscatter.
\end{abstract}

\section{Introduction}  \label{sec:intro}
This article describes sphere-to-sphere propagation of smooth
data for the ordinary 3-space dimensional wave equation 
\cite{Wilcox,Tokita,BFL,GHJ}.
As an application, consider evolving the wave equation \eqref{eq:3d}
on a finite computational domain with a spherical outer 
boundary of radius $r_1$, recording as a time-series 
the solution restricted to the boundary sphere, and 
--as a post-processing step-- recovering what the solution 
reaching $r_2 > r_1$ would be. In this application neither
the computational domain nor the final time need be extended.
For many applications, the post-processing step would be faster than 
evolving the wave equation on a commensurately larger spacetime
domain, whilst avoiding accumulation of phase or other 
errors typical of long-time integrations. We shall refer to 
the propagation of solution data from radius $r_1$ to radius
$r_2$ as {\em teleportation}, and in the limit 
$r_2 \rightarrow \infty$ as {\em asymptotic-waveform evaluation}.
As discussed in the concluding section, our methods likely 
extend to a certain class of hyperbolic equations.

To further elucidate the idea, consider evolution of the 
simple 1-space dimensional wave equation,
\begin{align} \label{eq:1p1}
-\partial_t^2 \Psi + \partial_x^2 \Psi = 0,
\end{align}
on $x_0 \le x \le x_1$, subject to Sommerfeld boundary 
conditions and with compactly supported initial data.
The general solution of d'Alembert is
\begin{align}\label{eq:1p1Sol}
\Psi(t,x) = F(x-t) + G(x+t),
\end{align}
where the "outgoing" wave is $F(x-t)$. As follows from 
the assumption of compact support and \eqref{eq:1p1Sol},
the wave value passing $x_1$ at time $t$ will reach 
$x_2 > x_1$ at time $t+(x_2-x_1)$, and
\begin{align}
\Psi(t + (x_2 - x_1), x_2) = F(x_1-t) =  \Psi(t, x_1).
\end{align}
This formula constitutes the simplest possible 
teleportation scheme, only accounting for the time-delay 
between the spacetime points $(t_1,x_1)$ and $(t_2,x_2)$.

Turn now to the 3-space dimensional wave equation,
\begin{align} \label{eq:3d}
(-\partial_t^2+\partial_x^2 + \partial_y^2 + 
  \partial_z^2)\psi = S(t,x,y,z),
\end{align}
subject to initial data and a source $S$ of compact support. We 
imagine a region of space enclosed by a spherical outer 
boundary of radius $r_1$ which need not be large. This boundary
sphere is also called the {\em extraction sphere}, and beyond 
it there are no scatters and both the initial data and source 
term vanish. Furthermore, on the boundary sphere we place 
exact radiation outer boundary 
conditions~\cite{AGH,GroteKeller1}.\footnote{An efficient 
simultaneous implementation of 
teleportation and radiation boundary conditions (RBC) would rely 
on common pole locations for both the teleportation and RBC kernels, 
thereby using the same ODEs for both recovery of the teleported 
signal and enforcement of the RBC. Reference~\cite{BFL} noted 
that, when achievable, the resulting teleportation kernels 
were of reduced accuracy. However, the preliminary study made
in Ref.~\cite{BFL} involved kernels for blackhole perturbations,
and future work should explore the issue for the wave equation.}
Let the Cartesian coordinates $\mathbf{x}=(x,y,z)$ be expressed as 
$\mathbf{x} = r\boldsymbol{\theta}$, with
$\boldsymbol{\theta}=(\sin\theta\cos\phi,\sin\theta\sin\phi,\cos\theta)$
the direction cosines associated to spherical polar 
coordinates.\footnote{Here with the convention that $\theta$ is the 
polar and $\phi$ the azimuthal angle.}
From the extraction-sphere data 
$r_1\psi(t,r_1\boldsymbol{\theta})$ we seek to recover the signal 
$r_2\psi(t+(r_2-r_1),r_2\boldsymbol{\theta})$ which would reach 
arbitrarily large distances $r_2$ (the extra radial factors 
account for the $1/r$ fall-off of $\psi$).

We adopt the following approach to far-field signal 
recovery: (i) derive an exact $r_1 \rightarrow r_2$ 
procedure/relationship in the spherical-harmonic-Laplace 
transform domain ("frequency domain"), and then (ii) approximate 
this exact relationship in a way 
which allows for a simple inversion of the Laplace 
transform~\cite{BFL}. Owing to the spatial spherical harmonic 
transform, the approach remains nonlocal in space. For the 
wave equation \eqref{eq:3d}, the relationship in the frequency 
domain (whether exact or approximate) involves a sum of simple 
poles; whence the corresponding time-domain procedure involves 
a history-dependent convolution based on a sum-of-exponentials 
kernel. For a given spherical-harmonic polar index $\ell$, the
{\em exact} sum-of-exponentials kernel involves precisely $\ell$
terms. As recently shown by Greengard, Hagstrom, and Jiang 
\cite{GHJ} and described further below, direct evaluation of the 
{\em exact} sum-of-poles frequency-domain kernel is 
problematic for large $\ell$. Indeed, in this case the 
complex residues in the sum vary in modulus over many orders 
of magnitude, and the sum is plagued by cancellation errors. 
Ref.~\cite{GHJ} has also described a method for evaluating 
the {\em time-domain} kernel via stable recursion relations, 
assuming that the pole locations (Bessel-MacDonald zeros) have 
been precomputed.

Using the technique of {\em kernel compression} 
\cite{AGH,Jiang,XuJiang,BEHL}, this paper considers
approximation of the exact frequency domain kernel.
Our approximations are based on an integral expression 
for the frequency domain kernel (cf.~Eq.~\eqref{eq:flatspacePhi})
which affords well-conditioned evaluation of the exact frequency 
domain kernel at imaginary Laplace frequencies, but by itself is not
useful for the time-domain. With expression \eqref{eq:flatspacePhi}
and a given integer $d \ll \ell$, we may construct an accurate 
rational approximation of the exact kernel which is itself a sum 
of $d$ simple poles, and therefore also determines a 
sum-of-exponentials kernel in the time-domain (now with fewer terms).
While offering no analytical proof,
we empirically demonstrate that kernel evaluation based on the
$d$-term approximate expansion is more accurate
than evaluation based on the $\ell$-term exact expansion.
Moreover, we examine scaling relations between the
parameters which determine a compressed kernel.
We also give a more detailed
derivation of Greengard, Hagstrom, and Jiang's large-$\ell$
asymptotic result \cite{GHJ} for the exact kernel residues.
Our approach was recently developed in
Ref.~\cite{BFL} for cases where closed-form kernel
expressions were unavailable, and the concluding section
remarks on applicability of the technique beyond the
ordinary wave equation case.

\section{Teleportation kernels}
This section reviews the origin and structure of teleportation
kernels \cite{Wilcox,Tokita,BFL,GHJ}, mostly following the 
presentation from \cite{BFL}. In the Laplace frequency $s$-domain, 
each such kernel is a finite sum of simple poles in the complex 
$s$-plane. Assembling various results from \cite{Olver,Olver2} and 
Abramowitz and Stegun's compendium \cite{AS} (hereafter AS), this 
section also considers the residues of a teleportation kernel, 
in particular deriving a large-$\ell$ asymptotic expansion. 
Our analysis in section \ref{subsec:asymp} is an elaboration of
results given in Ref.~\cite{GHJ}. 

\subsection{Derivation}
Assuming that $S=0$ (the source plays no direct role in what follows),
we start by expanding the solutions to Eq.~\eqref{eq:3d} as
\begin{equation} \label{eq:3dsol1}
\psi(t,x,y,z)= \frac{1}{r} \sum_{\ell = 0}^{\infty} \sum_{m=-\ell}^{\ell} 
\Psi_{\ell m}(t,r) Y_{\ell m}(\theta,\phi),
\end{equation}
where the $Y_{\ell m}(\theta,\phi)$ are standard spherical harmonics,
i.e.~the eigenfunctions of the Laplace-Beltrami operator on the unit 
sphere. Substitution of \eqref{eq:3dsol1} 
into \eqref{eq:3d} determines that for each $(\ell,m)$ pair the 
time-domain {\em multipole} $\Psi_{\ell m}(t,r)$ obeys
the {\em radial wave equation},
\begin{align}\label{eq:radialwave}
\partial_t^2 \Psi_{\ell m} - \partial_r^2 \Psi_{\ell m} 
+ \frac{\ell \left( \ell + 1 \right)}{r^2} \Psi_{\ell m} = 0.
\end{align}
Introducing the Laplace transform,
\begin{align}
\widehat{\Psi}_{\ell m}(s,r) = \int_0^{\infty} e^{-st}
\Psi_{\ell m}(t,r)dt,
\end{align}
we transform Eq.~\eqref{eq:radialwave}, with the result
\begin{align}\label{eq:MBS}
\left[s^2 - \frac{d^2}{d r^2} 
+ \frac{\ell(\ell + 1)}{r^2}\right]\widehat{\Psi}_{\ell m} 
= s\Psi_{\ell m}(0,r) 
+ \Psidot_{\ell m}(0,r). 
\end{align}
Solution of this equation by the method of variation of parameters
requires solutions to the homogeneous equation; these can be expressed in 
terms of modified Bessel functions.

Consider initial data of compact support, chosen to vanish 
outside of $r=r_1-\delta$, with $0 < \delta \ll 1$. Then for 
$r \geq r_1$, the general solution to \eqref{eq:MBS} is the
{\em outgoing} one
\begin{equation}\label{eq:ellmult_fd} 
\widehat{\Psi}_{\ell m}^{\tt out}(s,r) 
= A_{\ell m}(s)s^\ell e^{-sr} W_\ell(sr),
\qquad
W_\ell(z) = \sum_{k=0}^\ell \frac{c_{\ell k}}{z^k}, 
\qquad 
c_{\ell k} = \frac{1}{2^k k!}\frac{(\ell+k)!}{(\ell-k)!},
\end{equation}
where $z=sr$, $A_{\ell m}(s)$ is independent of $r$, and
(cf.~AS {\bf 10.1.9} and {\bf 10.2.15})
\begin{equation}\label{eq:KfromW}
W_\ell(z) = 
\sqrt{\frac{2z}{\pi}}e^z 
K_{\ell+1/2}(z).
\end{equation}
Here $K_\nu(z)$ is a modified spherical Bessel 
function (MacDonald's function) of Bessel order 
$\nu \equiv \ell + 1/2$ \cite{Watson,AS};
$\nu$ has this meaning throughout.

The structure of the Laplace-domain solution \eqref{eq:ellmult_fd}
determines the following algebraic relationship between solution values:
\begin{align} \label{eq:fdkernel}
\mathrm{e}^{s(r_2 - r_1)}
\widehat{\Psi}_{\ell m}^{\tt out}(s,r_2) = 
\widehat{\Phi}_\ell(s,r_1,r_2) 
\widehat{\Psi}_{\ell m}^{\tt out}(s,r_1) + 
\widehat{\Psi}_{\ell m}^{\tt out}(s,r_1),
\end{align}
where we have defined the frequency-domain teleportation kernel
(cf.~Eq.~(31) of \cite{BFL})
\begin{align} \label{eq:fdkernel_def}
\widehat{\Phi}_\ell(s,r_1,r_2) = -1 + 
\frac{W_\ell(sr_2)}{W_\ell(sr_1)}.
\end{align}
Here the "minus $1$" factor ensures that $\widehat{\Phi}_\ell(s,r_1,r_2)$
decays for large $s$, and therefore has an inverse Laplace transform which
is a classical function. Intuitively, the "minus $1$" factor also 
isolates $\widehat{\Phi}_\ell(s,r_1,r_2)\widehat{\Psi}_{\ell m}^{\tt out}(s,r_1)$ 
as the piece of the propagating wave which is altered as it moves from 
$r_1$ to $r_2$. Indeed, when $\ell = 0$ the kernel is zero and 
\eqref{eq:fdkernel} becomes
\begin{align}
e^{s(r_2 - r_1)} 
\widehat{\Psi}_{00}^{\tt out}(s,r_2)= 
\widehat{\Psi}_{00}^{\tt out}(s,r_1).
\end{align}
This formula involves only the wave transit time, and
could have been obtained directly in the time-domain 
(see the introduction). From Eq.~\eqref{eq:fdkernel_def}
the asymptotic signal for any $\ell$ is determined by
\begin{align} \label{eq:fdkernel_def_asy}
\widehat{\Phi}_\ell(s,r_1,\infty) 
= \frac{1 - W_\ell(sr_1)}{W_\ell(sr_1)},
\end{align}
where we have used $W_\ell(\infty) = 1$. The 
kernel $\widehat{\Phi}_{\ell}(s,r_1,r_2)$ teleports a 
signal $\widehat{\Psi}_{\ell m}^\mathrm{out}(s,r_1)$ of 
frequency $s$ from an extraction sphere of radius 
$r_1$ to $r_2 \leq \infty$. Expression~\eqref{eq:fdkernel} 
can be implemented in frequency-domain solvers. 
We now turn to the time-domain case.

The next subsection shows that \eqref{eq:fdkernel_def} 
can be represented as a sum of simple poles,
\begin{align} \label{eq:fdkernel_poles}
\begin{split}
\widehat{\Phi}_\ell(s,r_1,r_2)
= \sum_{j=1}^\ell \frac{a_{\ell j}(r_1,r_2)}{s - b_{\ell j}/r_1}.
\end{split}
\end{align}
Therefore, by well-known properties of the Laplace transform,
the inverse transformation of Eq.~\eqref{eq:fdkernel} is
(dropping the "out" superscript)
\begin{equation} \label{eq:td_teleportation} 
\Psi_{\ell m}(t+(r_2-r_1),r_2) = 
\int_0^t \Phi_\ell(t-t',r_1,r_2)\Psi_{\ell m}(t',r_1)dt' 
+ \Psi_{\ell m}(t,r_1),
\end{equation}
where the time-domain teleportation kernel is a sum of exponentials,
\begin{equation}\label{eq:tdPhi}
\Phi_\ell(t,r_1,r_2) 
= \sum_{k=1}^\ell a_{\ell k}(r_1,r_2)
e^{b_{\ell k} t/r_1}.
\end{equation}
The formulas~(\ref{eq:fdkernel_poles},\ref{eq:tdPhi}) express 
signal teleportation as a time-domain 
convolution, and for low-$\ell$ they are numerically useful.
However, as shown in \cite{GHJ} and reviewed below, for high $\ell$ 
they become exponentially ill-conditioned. Section \ref{sec:compress} 
discusses {\em kernel compression}, which involves approximation of
\eqref{eq:tdPhi} by a sum\footnote{Ref.~\cite{BFL} used $\Xi^E_\ell$,
$\gamma^E_{\ell,k}$, and $\beta^E_{\ell,k}$ for these quantities, 
where $E$ stands for "evaluation". In \cite{BFL} the same symbols 
{\em without} $E$ superscripts were used for similar quantities 
associated with radiation boundary conditions.}
\begin{equation}\label{eq:tdPhiE}
\Phi_\ell(t,r_1,r_2) \simeq 
\Xi_\ell(t,r_1,r_2) = \sum_{k=1}^d \gamma_{\ell k}(r_1,r_2)
e^{\beta_{\ell k}(r_1,r_2) t}
\end{equation}
of $d \leq \ell$ exponentials. 
Below we demonstrate empirically that
compression also alleviates the large-$\ell$ catastrophic cancellation
due to the exponential variation in size~\eqref{eq:mainresult} 
of the residues.

A convolution based on either \eqref{eq:tdPhi} or \eqref{eq:tdPhiE}
can be directly implemented in time-domain solvers, or carried out as a
post-processing step. The time series $\Psi_{\ell m}(t,r_1)$ must be
generated by a numerical solver, whereas the teleportation technique 
allows for reduction of the computational domain. That is, the 
signal $\Psi_{\ell m}(t+(r_2-r_1),r_2)$ {\em that would reach} 
$r_2$ is written explicitly in terms of data recorded on the 
extraction sphere. Appendix~\ref{app} provides an error
estimate for the teleportation of noisy numerical data. 
Efficiency gains in computing $r_2\psi(t,r_2\boldsymbol{\theta})$
by way of signal teleportation, as opposed to direct numerical simulation,
will depend on choices for $r_1$, $r_2$, the final simulation time, 
and the number of $\left(\ell,m\right)$ multipoles 
required to accurately resolve $r_1\psi(t,r_1\boldsymbol{\theta})$ 
with a truncated expansion~\eqref{eq:3dsol1}. These considerations 
are further discussed in Refs.~\cite{BFL,GHJ}.

\subsection{Frequency domain kernel as a sum of simple poles}
\begin{figure}[htp]
\begin{minipage}{0.4\textwidth}
\begin{center}
\includegraphics[width=10cm]{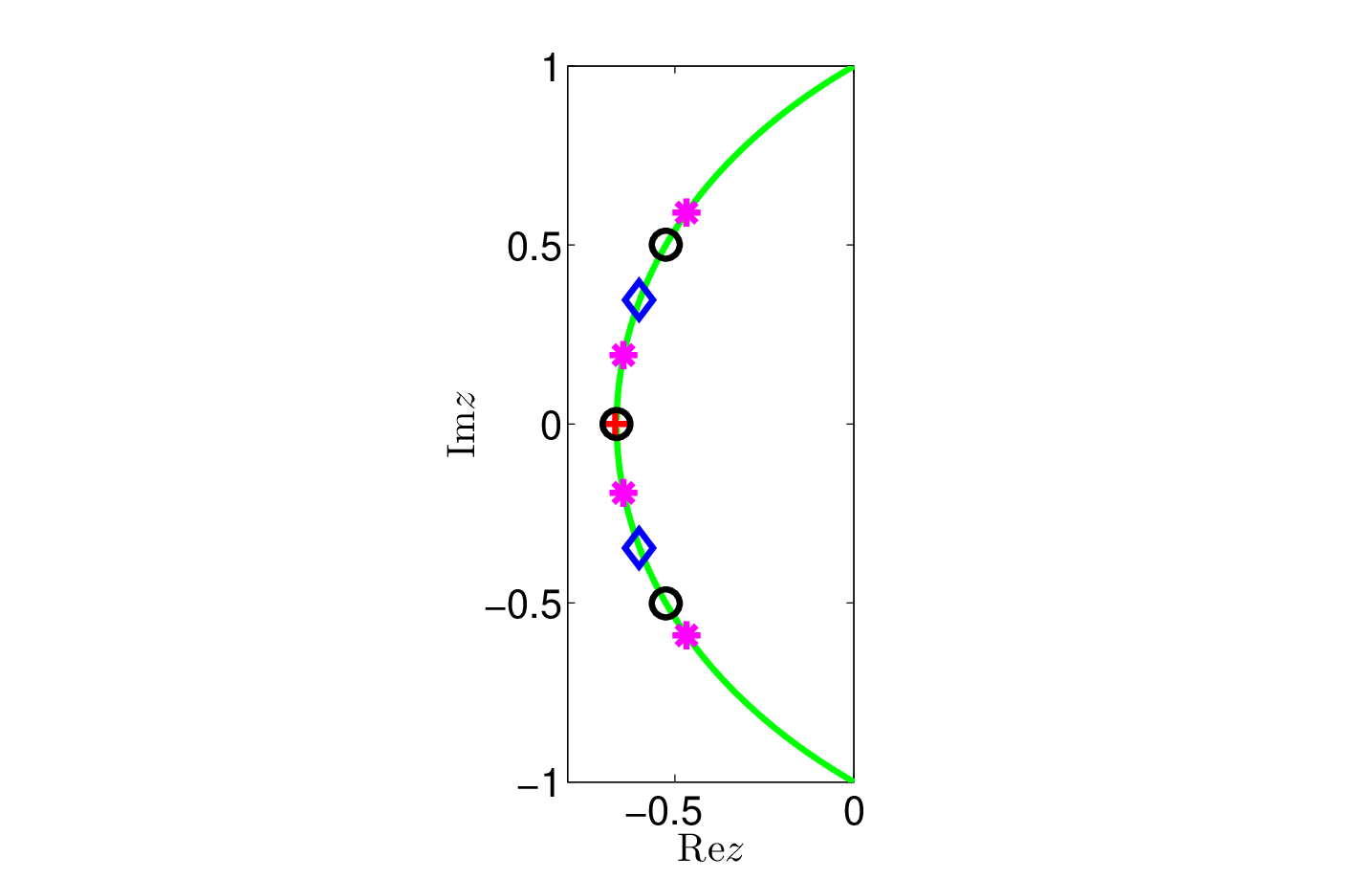}
\end{center}
\end{minipage}
\hspace{5mm}
\begin{minipage}{0.4\textwidth}
\begin{center}
{\large\begin{align*}
{\color{red}+} \displaystyle
z^1W_{1}(z) & = z
\left(1+\frac{1}{z}\right)
\\
{\color{blue}\diamond} \displaystyle
z^2W_{2}(z) & = z^2
\left(1+\frac{3}{z}+\frac{3}{z^2}\right)
\\
{\color{black}\circ} \displaystyle
z^3W_{3}(z) & = z^3
\left(1+\frac{6}{z}
             + \frac{15}{z^2} + \frac{15}{z^3}\right)
\\
{\color{magenta}*} \displaystyle
z^4 W_{4}(z) & = z^4
\left(1+\frac{10}{z} + \frac{45}{z^2}
+ \frac{105}{z^3} + \frac{105}{z^4}\right)
\end{align*}}
\end{center}
\end{minipage}
\caption{
Zeros $b_{\ell j}$ of $z^{\ell} W_{\ell}(z)$
{\em scaled} by order $\nu=\ell+1/2$.
Here we plot
$\{\tilde{b}_{\ell j}=b_{\ell j}/\nu : j=1,\dots,\ell\}$ for 
$1 \leq \ell \leq 4$, with the corresponding Bessel 
polynomials shown on the right. In the 
$\ell \rightarrow \infty$ limit the scaled zeros accumulate 
on a transcendental curve (green)~\cite{AGH,Olver}.
Evidently, the agreement holds even for the lowest
$\ell$, at least to the eye.}
\label{fig:PoleLocs}
\end{figure}

In terms of the Bessel-MacDonald zeros 
$\{b_{\ell j}: j = 1,\dots,\ell\}$ the function $W_\ell(z)$ has 
the form
\begin{align} \label{eq:W}
W_\ell(z) = z^{-\ell}\prod_{j=1}^\ell (z - b_{\ell j}) ,
\end{align}
which follows from \eqref{eq:KfromW} and shows that the 
roots of $W_\ell(z)$ are the same as those of $K_{\ell+1/2}(z)$.
The treatise by Watson shows that all $b_{\ell j}$ lie in 
the left-half plane $\{z : \mathrm{Re}z < 0\}$. 
Moreover, for even $\ell$ these roots come in $\ell/2$ conjugate 
pairs, while for odd $\ell$ we have $(\ell-1)/2$ conjugate 
pairs and precisely one negative real root.
In terms of the Bessel order $\nu$ the 
scaled roots $\tilde{b}_{\ell j} = b_{\ell j}/\nu$ are known to 
accumulate on a transcendental curve in the left-half 
plane~\cite{AGH,Olver}; see Fig.~\ref{fig:PoleLocs}.

From Eqs.~(\ref{eq:fdkernel_def},\ref{eq:W}) the kernel can be written as
\begin{align}
\begin{split}
\widehat{\Phi}_\ell(s,r_1,r_2)
& = -1 + \frac{\prod_{j=1}^\ell (s - b_{\ell j}/r_2)}{
               \prod_{j=1}^\ell (s - b_{\ell j}/r_1)}  
  = \sum_{j=1}^\ell \frac{a_{\ell j}(r_1,r_2)}{s - b_{\ell j}/r_1} .
\end{split}
\end{align}
Since the poles are simple, we compute the residues as
\footnote{The residues $a_{\ell j}(1,r)$ are precisely those 
considered in \cite{GHJ}.}
\begin{align} \label{eq:residues}
\begin{split}
a_{\ell j}(r_1,r_2) & = 
\lim_{s\rightarrow  b_{\ell j}/r_1} (s-b_{\ell j}/r_1) 
\widehat{\Phi}_\ell(s,r_1,r_2)
= r_1^{\ell-1}
  \frac{\prod_{k=1}^\ell (b_{\ell j}/r_1 - b_{\ell k}/r_2)}{
  \prod_{k=1, k\neq j}^\ell (b_{\ell j}- b_{\ell k})}.
\end{split}
\end{align}
The last expression in \eqref{eq:residues} implies that
$a_{\ell j}(r_1,r_2) = r_1^{-1}a_{\ell j}(1,r_2/r_1)$, from
which we infer the scaling relation
\begin{equation}\label{eq:scaling}
\widehat{\Phi}_\ell(s,r_1,r_2) = 
\widehat{\Phi}_\ell(sr_1,1,r_2/r_1).
\end{equation}

\subsection{Residue asymptotics}\label{subsec:asymp} 
Following Greengard, Hagstrom, and Jiang \cite{GHJ}, 
this subsection derives a large-$\ell$ asymptotic 
formula [\eqref{eq:mainresult} below] for the 
residues $a_{\ell j}(r_1,r_2)$. This asymptotic 
formula indicates that kernel evaluation based on 
the exact expressions 
(\ref{eq:fdkernel_poles},\ref{eq:td_teleportation})
is numerically inaccurate unless $\ell$ is small,
motivating our alternative procedure for
evaluation of frequency domain kernels described in 
Sec.~\ref{sec:evaluation}.
Since the time domain kernel
\eqref{eq:tdPhi} is real valued, the roots and residues
come in conjugate pairs. Whence it suffices to consider 
only those scaled roots $\tilde{b}_{\ell j} = b_{\ell j}/\nu$ 
obeying $-\pi \leq \arg \tilde{b}_{\ell j} < -\frac{1}{2}\pi$.
The asymptotic formula is
\begin{align}\label{eq:mainresult}
a_{\ell j}(r_1,r_2) \sim 
\frac{\mathrm{i}\tilde{b}_{\ell j}}{2r_1}
\left(\frac{r_2}{\pi r_1}\right)^{1/2}
\Big(
1+\tilde{b}_{\ell j}^2
\Big)^{-1/4}
\Big(
1+\tilde{b}_{\ell j}^2(r_2/r_1)^2
\Big)^{-1/4}
\frac{
e^{\nu\psi_{\ell j}(r_1,r_2)}
}{
(a_j)^{-1/4}
\mathrm{Ai}'(a_j) 
},\qquad \ell \rightarrow \infty .
\end{align}
In this expression 
\begin{equation}
\psi_{\ell j}(r_1,r_2) = (r_2/r_1 - 1)\tilde{b}_{\ell j} -
\log\frac{1+\sqrt{1+\tilde{b}_{\ell j}^2(r_2/r_1)^2}
}{
\mathrm{i} \tilde{b}_{\ell j}r_2/r_1}
+ \sqrt{1+\tilde{b}_{\ell j}^2(r_2/r_1)^2} \, ,
\end{equation}
and $a_j$ is the $j$th root of the Airy function $\mathrm{Ai}(y)$.
To achieve the correct correspondence between the 
$\tilde{b}_{\ell j}$ and the $a_j$, formula \eqref{eq:mainresult} 
assumes that the set $\{\tilde{b}_{\ell j} : j=1,\dots,\ell\}$ is 
ordered from the bottom of the third quadrant upwards, 
i.e.~for $\ell = 3$ the (unscaled) roots $b_{3j}$ are ordered as 
follows (of these we consider only $b_{31}$ and $b_{32}$).
\begin{verbatim}
     b31 = -1.8389e+00 - 1.7544e+00i
     b32 = -2.3222e+00 + 0.0000e+00i
     b33 = -1.8389e+00 + 1.7544e+00i 
\end{verbatim}
\begin{remark} \label{remark1}
At the expense of introducing a second expansion with its own
error, \eqref{eq:mainresult} may be further reduced through
the following large-$j$ formulae (see the appendix of
\cite{Olver2}):
\begin{equation}
a_j \sim
-(\threehalf \pi )^{2/3}(j-\forth)^{2/3},
\quad
\mathrm{Ai}'(a_j) \sim \frac{(-1)^{j-1}}{\sqrt{\pi}}
(\threehalf \pi)^{1/6}(j-\forth)^{1/6},
\quad
(a_j)^{-1/4}\mathrm{Ai}'(a_j) \sim
\frac{(-1)^{j-1}}{\sqrt{\pi}}
e^{-\mathrm{i}\pi/4}.
\end{equation}
\end{remark}

We now turn to the derivation of \eqref{eq:mainresult}.
Rather than \eqref{eq:residues}, the derivation starts with \cite{GHJ}
\begin{equation} \label{eq:afromscaled}
a_{\ell j}(r_1,r_2) =
-\frac{\mathrm{i}}{r_1}\sqrt{\frac{r_2}{r_1}}
e^{(r_2/r_1 - 1)\nu\tilde{b}_{\ell j}}
\frac{
H_{\nu}^{(1)}(\mathrm{i}\nu \tilde{b}_{\ell j} r_2/r_1)
}{
H_{\nu}^{(1)\prime}(\mathrm{i}\nu \tilde{b}_{\ell j})
}.
\end{equation}
Appendix \ref{sec:appendixA} shows how this formula stems
from \eqref{eq:residues}.
Since we consider those $\mathrm{i}b_{\ell j}$ in the fourth 
quadrant, the arguments of the Hankel functions in 
\eqref{eq:afromscaled} are certainly in the sector of validity 
for the following expansions (AS {\bf 9.3.37} and {\bf 9.3.45}):
\begin{subequations}\label{eq:olver}
\begin{align}
H_\nu^{(1)}(\nu w) & \sim 2 e^{-\mathrm{i}\pi/3}\nu^{-1/3}
\left(\frac{4\zeta}{1-w^2}\right)^{1/4}
\mathrm{Ai}(e^{2\pi\mathrm{i}/3}\nu^{2/3} \zeta)\\
H_\nu^{(1)\prime}(\nu w) & \sim \frac{4}{w} e^{4\mathrm{i}\pi/3}\nu^{-2/3}
\left(\frac{1-w^2}{4\zeta}\right)^{1/4}
\mathrm{Ai}'(e^{2\pi\mathrm{i}/3}\nu^{2/3}\zeta),
\end{align}
\end{subequations}
which hold for $|\arg w| \leq \pi - \epsilon$. 
In these expansions
(see AS {\bf 9.3.38}, but note that the $z$ and $\ln$ in AS are our 
$w$ and $\log$)
\begin{equation}\label{eq:magic}
\frac{2}{3}\zeta^{3/2} = 
\log\frac{1+\sqrt{1-w^2}}{w} - \sqrt{1-w^2},
\end{equation}
defines $\zeta$ implicitly as function of $w$.

We first consider evaluation of (\ref{eq:olver}b) with 
$w = \mathrm{i}\tilde{b}_{\ell j}$. As
$\nu \rightarrow \infty$ the scaled roots 
$\tilde{b}_{\ell j}$ obey 
\cite{Olver}
\begin{equation}
\mathrm{i}\tilde{b}_{\ell j}\sim 
w(\zeta_j) + O(\nu^{-1}),\qquad
\zeta_j = e^{-2\pi\mathrm{i}/3}\nu^{-2/3} a_j,
\end{equation}
uniformly in $j$.
Therefore, the desired expression for (\ref{eq:olver}b) is
\begin{equation}\label{eq:HankelPrime_asym}
H_\nu^{(1)\prime}(\nu \mathrm{i}\tilde{b}_{\ell j})
\sim 
-2\sqrt{2}\nu^{-1/2}(\tilde{b}_{\ell j})^{-1}
(1+\tilde{b}_{\ell j}^2)^{1/4}(a_j)^{-1/4}
\mathrm{Ai}'(a_j).
\end{equation}
In canceling terms to reach this expression, we have 
paid due attention to the branch associated with the fourth root. 

Next, we turn to (\ref{eq:olver}a) with 
$w = \mathrm{i}b_{\ell j}r_2/r_1$,
\begin{align}
H_\nu^{(1)}(\nu \mathrm{i}b_{\ell j}r_2/r_1) 
& \sim 2 e^{-\mathrm{i}\pi/3}\nu^{-1/3}
\left(\frac{4\zeta_r
}{
1+\tilde{b}_{\ell j}^2(r_2/r_1)^2
}
\right)^{1/4}
\mathrm{Ai}(e^{2\pi\mathrm{i}/3}\nu^{2/3}\zeta_r),
\end{align}
where $\zeta_r = \zeta(\mathrm{i}b_{\ell j}r_2/r_1)$, 
that is $\zeta_r$ is stems from \eqref{eq:magic} 
with $w = \mathrm{i}b_{\ell j}r_2/r_1$. 
The analysis given in \cite{Olver2} 
(pages 355-356) shows that $\frac{1}{3}\pi < \arg(\zeta_r)
< \pi$, provided $w$ is in fourth quadrant. 
The large argument asymptotics 
$\mathrm{Ai}(y)\sim \half\pi^{-1/2} y^{-1/4}
e^{-\frac{2}{3}y^{3/2}}$
of the Airy function determine that
\begin{align}
\mathrm{Ai}(e^{2\pi\mathrm{i}/3}\nu^{2/3} \zeta)
& \sim 
\half \pi^{-1/2} \nu^{-1/6} e^{\pi\mathrm{i}/3}\zeta^{-1/4}
e^{-\frac{2}{3}\nu \zeta^{3/2}}.
\end{align}
To reach this equation from the Airy expansion, 
we have assumed $\frac{1}{3}\pi < \arg \zeta < \pi$
and $-\pi < \arg(e^{2\pi\mathrm{i}/3}\zeta) < -\frac{1}{3}\pi$.
Recall that $\zeta_r$ obeys these inequalities.
With the last result and further consideration of the
fourth-root branch, we get
\begin{align}
H_\nu^{(1)}(\nu \mathrm{i}b_{\ell j}r_2/r_1)
& \sim \sqrt{2}\pi^{-1/2}\nu^{-1/2}
\Big(1+\tilde{b}_{\ell j}^2(r_2/r_1)^2\Big)^{-1/4}
e^{-\frac{2}{3}\nu \zeta_r^{3/2}},
\end{align}
and finally
\begin{align}
\frac{H_\nu^{(1)}(\nu \mathrm{i}b_{\ell j}r_2/r_1)
}{
H_\nu^{(1)\prime}(\nu \mathrm{i}\tilde{b}_{\ell j})
}
\sim 
-\frac{\tilde{b}_{\ell j}}{2\sqrt{\pi}}
\Big(
1+\tilde{b}_{\ell j}^2
\Big)^{-1/4}
\Big(
1+\tilde{b}_{\ell j}^2 (r_2/r_1)^2
\Big)^{-1/4} 
\frac{
e^{-\frac{2}{3}\nu \zeta_r^{3/2}}
}{
(a_j)^{-1/4}
\mathrm{Ai}'(a_j)
}.
\end{align}
This equation and \eqref{eq:afromscaled} yield 
\eqref{eq:mainresult}.

\section{Numerical approximation of 
teleportation kernels}\label{sec:numKernels}

This section treats the numerical approximation of 
frequency-domain teleportation kernels $\widehat{\Phi}_\ell(s,r_1,r_2)$.
Construction of our approximations requires that we are able to 
accurately evaluate $\widehat{\Phi}_\ell(\mathrm{i}y,r_1,r_2)$ for any 
$y\in\mathbb{R}$. Due to the variation 
in size of the residues for high $\ell$,
Eq.~\eqref{eq:fdkernel_poles} does not offer a viable means
for such evaluation (even given the ability 
to compute the poles and residues), 
unless $\ell$ is low. However, numerical evaluation 
based on Eq.~\eqref{eq:flatspacePhi} (below) is well 
conditioned. Use of Eq.~\eqref{eq:flatspacePhi} amounts 
to a cumbersome offline step, through which we construct an 
accurate rational approximation $\widehat{\Xi}_\ell(s,r_1,r_2)$ 
to $\widehat{\Phi}_\ell(s,r_1,r_2)$ along the inversion contour.
The approximate kernel $\widehat{\Xi}_\ell(s,r_1,r_2)$ 
typically has $d \ll \ell$ poles for large $\ell$, 
and is therefore a "compression" of the exact $\ell$-pole 
kernel. Quite remarkably, (approximate) evaluation of 
$\widehat{\Phi}_\ell(\mathrm{i}y,r_1,r_2)$ based on the $d$-pole
sum $\widehat{\Xi}_\ell(\mathrm{i}y,r_1,r_2)$ is much better
conditioned than evaluation based on the exact 
$\ell$-pole sum \eqref{eq:fdkernel_poles}.

\subsection{Profile evaluation}\label{sec:evaluation}
We now describe our alternate approach for evaluation of the
profiles $\mathrm{Re}\widehat{\Phi}_\ell(\mathrm{i}y,r_1,r_2)$ and 
$\mathrm{Im}\widehat{\Phi}_\ell(\mathrm{i}y,r_1,r_2)$ for $y\in\mathbb{R}$.
As mentioned, Ref.~\cite{GHJ} has described stable evaluation 
of the {\em time-domain} kernel $\Phi_\ell(t,r_1,r_2)$.
Our approach is based on the following expression for a teleportation
kernel~\cite{BFL}:
\begin{equation}\label{eq:flatspacePhi}
\widehat{\Phi}_\ell(s,r_1,r_2)
         = -1 + \underbrace{\exp\left[\int_{r_1}^{r_2}
           \frac{\widehat{\Omega}_\ell(s,\eta)}{\eta}d\eta\right]}_{
           W_\ell(sr_2)/W_\ell(sr_1)},
\end{equation}
where we have introduced an auxiliary function\footnote{The kernel
$\widehat{\Omega}_\ell(s,r)$ arises when deriving exact 
outgoing (i.e.~non-reflecting) boundary conditions. This interesting 
relationship  expresses a teleportation kernel as a weighed 
integral over boundary kernels.}
\begin{equation} \label{eq:Omegahat}
\widehat{\Omega}_\ell(s,r)
\equiv sr \frac{W'_\ell(sr)}{W_\ell(sr)}
= \sum_{k=1}^\ell \frac{b_{\ell k}/r}{s-b_{\ell k}/r},
\end{equation}
with the prime indicating differentiation in argument. 
With Steed's algorithm \cite{Thompson} the kernel 
$\widehat{\Omega}_\ell(s,r)$ is accurately 
computed via the known continued fraction expansion
\cite{Thompson}
\begin{equation}\label{eq:continuedfraction}
z \frac{W_\ell'(z)}{W_\ell(z)} = -\frac{\ell(\ell+1)}{2(z+1)+}
\frac{(\ell-1)(\ell+2)}{2(z+2)+} \cdots
\frac{2(2\ell-1)}{2(z+\ell-1)+} \frac{2\ell}{2(z+\ell)}.
\end{equation}
This formula follows from recurrence relations obeyed by 
MacDonald functions \cite{Watson}. Given the ability to compute
$\widehat{\Omega}_\ell(\mathrm{i}y,r)$, computation of 
\eqref{eq:flatspacePhi} can be carried out using numerical 
quadrature. Due to the structure of
$\mathrm{Re}\widehat{\Omega}_\ell(\mathrm{i}y,r)$ and
$\mathrm{Im}\widehat{\Omega}_\ell(\mathrm{i}y,r)$, the
numerical integrations involve no cancellation errors (i.e.~the sums
involve only positive or negative values at each fixed
$y_j$ grid point). We have found the 
representation~\eqref{eq:flatspacePhi} a useful tool 
for evaluation of high-$\ell$ teleportation kernels when evaluation based on
the sum-of-poles representation~\eqref{eq:fdkernel_poles} 
is inaccurate. Indeed, using quadruple precision arithmetic,
we are typically able to evaluate the profiles
$\mathrm{Re}\widehat{\Phi}_\ell(\mathrm{i}y,r_1,r_2)$ and
$\mathrm{Im}\widehat{\Phi}_\ell(\mathrm{i}y,r_1,r_2)$ with 
relative errors well below double precision accuracy.

\subsection{Compression}\label{sec:compress}
{\em Kernel compression} involves 
approximation of $\widehat{\Phi}_\ell(s,r_1,r_2)$ as a sum of 
(typically far) fewer poles. More precisely, given a prescribed error
tolerance $\varepsilon$, consider the approximation\footnote{
We correct several typos in Ref.~\cite{BFL}. In Eq.~(37) of 
that reference, the second summation sign $\Sigma$ should be 
a $\mathrm{sup}$ (supremum). Also in line 3 of {\bf Algorithm 4}, 
each $\widehat{\omega}_2$ should be $\widehat{\omega}_\ell$.
}

\begin{equation}\label{eq:compressedFLTkern}
\widehat{\Xi}_\ell(s,r_1,r_2) =
\sum_{n=1}^d \frac{\gamma_{\ell n}(r_1,r_2)}{s-\beta_{\ell n}(r_1,r_2)},
\qquad
\sup_{s\in \mathrm{i}\mathbb{R}}
\left|\frac{\widehat{\Phi}_\ell(s,r_1,r_2)-\widehat{\Xi}_\ell(s,r_1,r_2)}{
\widehat{\Phi}_\ell(s,r_1,r_2)}\right| < \varepsilon.
\end{equation}
In fact, due to scaling relation \eqref{eq:scaling} it suffices
to consider only approximations for $1\rightarrow r_2/r_1$ 
teleportation.\footnote{Indeed, other scenarios are 
determined by the rules $\gamma_{\ell n}(r_1,r_2) = 
r_1^{-1}\gamma_{\ell n}(1,r_2/r_1)$
and $\beta_{\ell n}(r_1,r_2) = \beta_{\ell n}(r_1) 
= r_1^{-1}\beta_{\ell n}(1)$.}
However, in practice we construct compressions of generic $r_1 
\rightarrow r_2$ kernels.

We now provide an error estimate associated with performing 
the convolution \eqref{eq:td_teleportation} with the 
compressed kernel $\Xi_\ell(t,r_1,r_2)$ given by 
\eqref{eq:tdPhiE} in place of the exact kernel 
$\Phi_\ell(t,r_1,r_2)$ given by Eq.~\eqref{eq:tdPhi}. 
From the Parseval and Fourier convolution theorems, the relative 
convolution error stemming from a compressed 
kernel is~\cite{AGH,BFL}
\begin{align}\label{eq:basic_estimate_rel}
\begin{split}
\big\|\Xi_\ell(\cdot,r_1,r_2)*\Psi_{\ell m}(\cdot,r_1)
& -\Phi_\ell(\cdot,r_1,r_2)*\Psi_{\ell m}(\cdot,r_1)\big\|_{L_2(0,\infty)} \\
& \leq \sup_{s\in\mathrm{i}\mathbb{R}}
\frac{|\widehat{\Xi}_\ell(s,r_1,r_2)-\widehat{\Phi}_\ell(s,r_1,r_2)|
}{
|\widehat{\Phi}_\ell(s,r_1,r_2)|} \times
\big\|\Phi_\ell(\cdot,r_1,r_2)*\Psi_{\ell m}(\cdot,r_1)\big\|_{L_2(0,\infty)}.
\end{split}
\end{align}
Assuming a rational approximation which achieves 
\eqref{eq:compressedFLTkern}, combination of 
\eqref{eq:compressedFLTkern} and \eqref{eq:basic_estimate_rel}
shows that the tolerance $\varepsilon$ is a long-time
bound on the relative convolution error.

We briefly describe our construction of compressed teleportation kernels.
Provided that $r_2-r_1$ is not too large, we use Algorithm 
\ref{alg:compressTLB} to produce a compressed kernel 
$\widehat{\Xi}_\ell(s,r_1,r_2)$. This is AGH compression, 
as described in \cite{AGH,Jiang,XuJiang,BEHL}. All of the
kernels used in this article have been generated via a
quadruple-precision implementation of this algorithm.
\begin{algorithm}{
{\sc Computation of a compressed teleportation kernel.}
\label{alg:compressTLB}
\vskip 5pt
\noindent
INPUT: $\ell$, $r_1$, $r_2$, $\varepsilon$, 
$N_\mathrm{C}$ (number of composite subintervals), 
$d$ (desired number of poles, possibly updated)
\\
OUTPUT: $\{\beta_{\ell n}(r_1,r_2),\gamma_{\ell n}(r_1,r_2)\}_{n=1}^d$}
\vskip 5pt
\begin{algorithmic}[1]
\State Choose an approximation window
       $[-y_\mathrm{max},y_\mathrm{max}]$
       on the $\sigma = \mathrm{i}y$ imaginary axis.
\State Partition $[-y_\mathrm{max},y_\mathrm{max}]$ to
       form a $y$-grid $\{y_j\}_{j=1}^J$, 
       typically with mesh refinement at the origin.
\State
      Numerically evaluate the profiles
      $\mathrm{Re}\widehat{\Phi}_\ell(\mathrm{i}y_j,r_1,r_2)$
      and $\mathrm{Im}\widehat{\Phi}_\ell(\mathrm{i}y_j,r_1,r_2)$
      on the $y$-grid via approximation of \eqref{eq:flatspacePhi}
      using $N_C$-composite Gauss-Kronrod quadrature.
\State
      Compute the numbers $\{\beta_{\ell n}(r_1,r_2),
      \gamma_{\ell n}(r_1,r_2)\}_{n=1}^d$
      by AGH compression; see Ref.~\cite{AGH}. 
      The idea is to solve
      \begin{equation} \label{eq:LS_AGH}
      \min_{\{\beta_{\ell n}(r_1,r_2),\gamma_{\ell n}(r_1,r_2)\}_{n=1}^d}
      \hspace{6pt}
      \sum_{j=1}^J\mu_j\Big|
      \widehat{\Phi}_\ell(\mathrm{i}y_j,r_1,r_2) - 
      \sum_{n=1}^d \gamma_{\ell n}(r_1,r_2)\big/
      \big(\mathrm{i}y_j - \beta_{\ell n}(r_1,r_2)\big)
      \Big|^2,
      \end{equation}
      where $\mu_j$ are quadrature weights.
\State
      Using 
      $\{\beta_{\ell n}(r_1,r_2),\gamma_{\ell n}(r_1,r_2)\}_{n=1}^d$,
      verify \eqref{eq:compressedFLTkern} for the chosen $\varepsilon$
      (typically on a different and finer $y$-grid). 
      If not verified, repeat last two steps with $d \leftarrow d+1$.
\end{algorithmic}
\end{algorithm}

For reasons discussed in \cite{BFL}, when $r_2 \gg r_1$ Algorithm
\ref{alg:compressTLB} becomes impractical. We therefore describe a
modified procedure, assuming for simplicity that $r_1 = 1$ and
$r_2 = 10^{P+1}$. We first compute the compressed 
kernel $\widehat{\Xi}^{(0)}_\ell(s,1,10)
\equiv \widehat{\Xi}_\ell(s,1,10)$ using 
Algorithm \ref{alg:compressTLB}.
As suggested by the scaling relation 
\eqref{eq:scaling}, teleportation over 
decades $[10^p,10^{p+1}]$ is then defined through
\begin{equation}
\widehat{\Xi}^{(p)}_\ell(s,10^{p},10^{p+1})
= \widehat{\Xi}^{(0)}_\ell(s 10^p,1,10) 
= \sum_{n=1}^d \frac{ 10^{-p}\gamma_{\ell n}(1,10)
}{ 
s-10^{-p}\beta_{\ell n}(1,10)}.
\end{equation}
Finally, we perform the evaluations more cheaply in step 3
of Algorithm \ref{alg:compressTLB} using a different formula,
\begin{equation}
\widehat{\Phi}_\ell(\mathrm{i}y_j,1,10^{P+1})
\simeq -1 + \prod_{p=0}^P 
\big[1+\widehat{\Xi}^{(p)}_\ell(\mathrm{i}y_j,10^{p},10^{p+1})].
\end{equation}
Here a collection of compressed kernels is combined to generate
profiles from which one new compressed kernel is obtained.

\subsection{Error estimate for 
pole-sum approximations}\label{sec:pole_approx}

We have used Algorithm~\ref{alg:compressTLB} (or 
its elaboration discussed in the last subsection) to 
construct compressed kernels. In view of the error estimate 
\eqref{eq:basic_estimate_rel}, we wish to know (or bound) 
the maximum pointwise error associated with a compressed 
kernel. To estimate this error in practice, we have 
resorted to numerical comparison between the compressed 
kernel and the "truth" kernel; see 
Subsection~\ref{sec:num_compression}. However, here we 
describe an a priori estimate for the relevant pointwise error. 

Consider a sum $f(z)$ of $n$ simple poles. The 
function $f(z)$ may represent one of our frequency domain 
kernels, although here we use $z$ in place of $s$ as the 
independent variable. Provided that evaluation of $f(z)$ is 
restricted to $z$-values which are sufficiently separated
from the pole locations $\{z_i\}_{i=1}^n$, a sum $g(z)$ of $d < n$ 
simple poles may accurately approximate $f(z)$. Reference 
\cite{AGH} analyzes this issue, giving both explicit 
constructions and error bounds. Such constructions are not 
used in Algorithm~\ref{alg:compressTLB}. Nevertheless, they 
demonstrate that good approximations 
exist and may shed light on observed numerical behavior 
(cf.~Subsection \ref{sec:num_scalings}).

\begin{lemma} \label{lem:pole_approx}
Suppose $\{q_j,z_j\}_{j=1}^n$ are $n$ complex pairs defining 
the function
$f(z) = \sum_{j=1}^n q_j/(z-z_j)$,
with all locations $z_j$ contained in the union of $p$ 
disks $D_1, \dots, D_p$. The disk $D_k$ has radius $r_k$ 
and is centered at $c_k$. There exists an approximation 
$g(z)$, itself a sum of $d = m\cdot p$ simple poles, 
which obeys the estimate
\begin{equation} \label{eq:pointwise_pole_err}
| f(z) - g(z) | \leq 
\frac{K (a^2+1)}{(a^m-1)(a-1)^2} 
\Big| \sum_{j=1}^n \frac{|q_j|}{z-z_j}\Big|,\qquad
z \in 
\mathcal{U}_a \equiv \big\{ z \big|
\mathrm{Re}(z - c_k)
\geq a r_k > r_k, 1 \leq k \leq p \big\},
\end{equation}
where the constant $K$ is independent of $a$.
\end{lemma}

For disks in the left-half plane we choose the largest $a > 1$ such that
$\mathcal{U}_a$ contains the inversion contour (i.e.~the imaginary axis).
Provided $n$ is sufficiently large, the estimate \eqref{eq:pointwise_pole_err} 
indicates that the error decays exponentially with the number $m\cdot p$ 
of poles. For an index-$\ell$ radiation boundary kernel, the residues
and pole locations are $q_j = b_{\ell j} = z_j$, assuming a 
unit-sphere physical boundary. Since $b_{\ell j} \sim \ell$ for large-$\ell$, 
these kernels correspond to a bound \eqref{eq:pointwise_pole_err} which 
scales mildly with $\ell$. For an index-$\ell$ teleportation kernel 
(with $r_1 = 1$) the locations are again $z_j =  b_{\ell j}$; however,
now the residues $q_j = a_{\ell j}$ scale exponentially with $\ell$, 
maring the bound \eqref{eq:pointwise_pole_err}.
Nevertheless, Section~\ref{sec:num} demonstrates that 
Algorithm~\ref{alg:compressTLB} yields accurate large-$\ell$ 
approximate kernels.

The discussion above serves to raise the following points. 
(i) Despite the fact that the kernels 
for radiation boundary conditions and teleportation share pole locations, 
teleportation kernels should be harder to approximate. (ii) Regardless, 
for fixed $\ell$, approximations with exponential accuracy exist in 
principle. (iii) For both types of kernels the set of pole locations 
can be covered by $p \propto \log \ell$ disks (Ref.~\cite{AGH} uses
this $p$ for RBC kernels). (iv) In the bound 
\eqref{eq:pointwise_pole_err} the exponential growth of the
residues $q_j = a_{\ell j}$ with $\ell$ suggests that to 
achieve an approximation with a fixed error tolerance the
number $d$ of poles needs to scale linearly with $\ell$. 
Subsection~\ref{sec:num_scalings} presents numerical evidence
supporting this heuristic observation. A problem 
whose solution might address the conditioning of 
our approximations would be to bound the largest and smallest 
residues (in modulus) for a given $g(z)$, and further to 
understand how these bounds scale with $\ell$.
Subsection \ref{sec:num_scalings} empirically addresses some of 
these issues.

\section{Example compressed kernels, empirical 
scalings, and numerical experiments} \label{sec:num}

Our first subsection 
describes approximation of large-$\ell$ frequency domain kernels. 
Our goal is the construction of compressed kernels, each with 
poles and residues that, when stored in double precision scientific 
format, yield accurate evaluation of the kernel along the inversion 
contour. Our second subsection summarizes empirical
scaling relations between $\ell$, $d$, $r_2$, $r_1$,  
$\max_j |\gamma_{\ell j}|$, and $\varepsilon$. 
In the last subsection two experiments describe signal teleportation. 
These experiments evolve time-domain multipole $\Psi_{\ell m}(t,r)$ 
which obey \eqref{eq:radialwave}. Working with the one-dimensional radial 
wave equation allows us to more easily focus on the errors associated
with the kernels as opposed to computational grid discretization error. 
Teleportation schemes implemented within three-dimensional time-domain 
wave equation solvers would require spherical harmonic transformation
to compute each multipole.

\subsection{Example compressions of large-$\ell$ kernels} \label{sec:num_compression}
\begin{figure}[htp]
\begin{center}
\subfigure[Pan out.]{
\includegraphics[width=7.5cm]{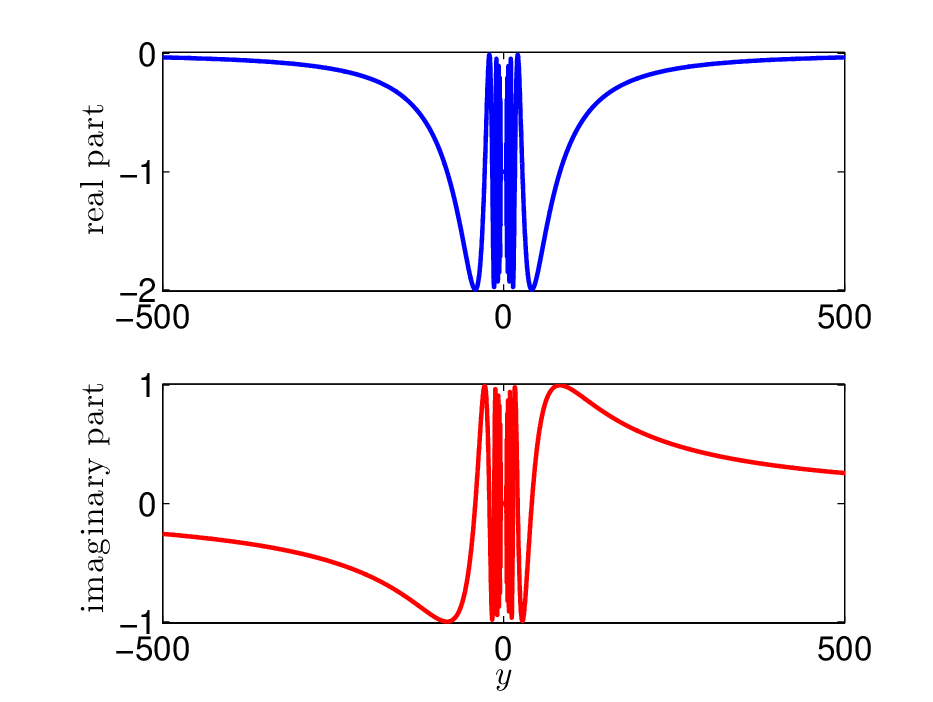}
}
\subfigure[Zoom in.]{
\includegraphics[width=7.5cm]{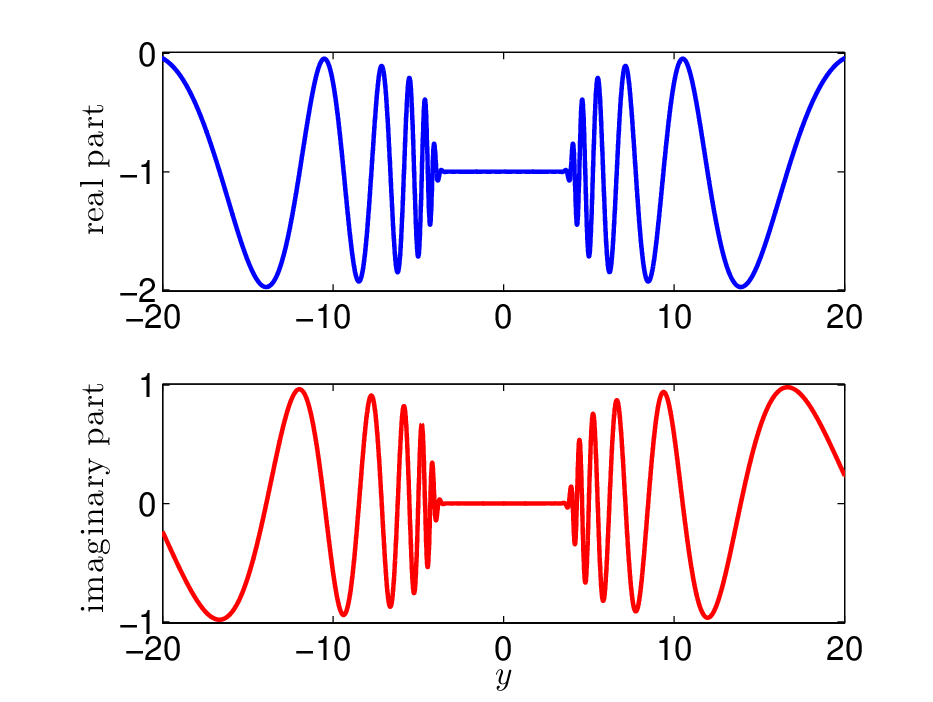}
}
\end{center}
\caption{Teleportation kernel $\widehat{\Phi}_{64}(\mathrm{i}y,15,240)$.
The upper plots depict Re$\widehat{\Phi}_{64}(\mathrm{i}y,15,240)$
and the bottom plots Im$\widehat{\Phi}_{64}(\mathrm{i}y,15,240)$.}
\label{fig:ell64kernel}
\end{figure}

We first consider $\ell=64$, with teleportation from either
$r_1=15$, $30$, $60$, or $120$ to $r_2=240$.
Due to the scaling relation \eqref{eq:scaling} 
these cases are similar to teleportation from $r_1=1$ to
$r_2 =16,8,4,2$. Figure \ref{fig:ell64kernel} shows the real 
and imaginary profiles associated with the kernel 
$\widehat{\Phi}_{64}(s,15,240)$ along the inversion contour.
Based on 
\eqref{eq:mainresult}, we expect that the exact residues 
$a_{64,j}(15,240)$ corresponding to the kernel 
$\widehat{\Phi}_{64}(s,15,240)$ vary in modulus over 15 
orders of magnitude, with the largest residue 
$\gtrsim 10^{15}$ in modulus. Therefore, one expects that a
table of the poles and residues specifying 
$\widehat{\Phi}_{64}(s,15,240)$ would, if formatted in
double precision, yield few digits of accuracy. 
Nevertheless, perhaps somewhat paradoxically, we will demonstrate
that the kernel $\widehat{\Phi}_{64}(\mathrm{i}y,15,240)$ can be
uniformly approximated for $y\in\mathbb{R}$ by a sum of 
{\em fewer poles}. 

Table \ref{tab:ell64} summarizes our best compressions for 
all four choices of $r_1$, and we find that between 
$30$ and $34$ poles are sufficient to
achieve $\varepsilon \lesssim$ {\tt 1.0e-10}.
To estimate the error tolerance
$\varepsilon$ in \eqref{eq:compressedFLTkern}
each approximate kernel is compared
against the "truth" kernel (i.e.~numerical evaluation of
$\widehat{\Phi}_{64}(\mathrm{i}y,r_1,240)$ with
\eqref{eq:flatspacePhi}) relative to a 
dense reference grid with logarithmic refinement.
In estimating $\varepsilon$, we also change the parameters 
for numerical quadrature used in evaluation of the 
"truth" kernel, thereby avoiding 
systematic errors. Reference \cite{BFL} has also examined
compressed kernels approximating $\widehat{\Phi}_{64}(s,15,240)$,
in particular plotting the pole locations for $d=20$,
28, and 36 compressed kernels. Figure 3 of that reference 
compares the pole locations of compressed teleportation and 
RBC kernels. For both types of compressed kernel as $d$ 
increases (corresponding to smaller $\varepsilon$),
more of the $d$ pole locations "lock on" to the transcendental 
curve shown in Fig.~\ref{fig:PoleLocs}.
\begin{table*}[h!]
\begin{tabular}{c || c c c c c c}
$r_1$ & $d$ & $\varepsilon$ & $\min_j |\gamma_{64,j}(r_1,240)|$ & $\max_j |\gamma_{64,j}(r_1,240)|$ &
$\min_j |a_{64,j}(r_1,240)|$ & $\max_j |a_{64,j}(r_1,240)|$\\
\hline
$15$   & $34$ & {\tt 1.0e-10} & {\tt 8.5972e-03} & {\tt 1.3898e+06} & {\tt 1.9086e+01} &  {\tt 4.3234e+15} \\
$30$   & $32$ & {\tt 1.0e-11} & {\tt 1.1253e-02} & {\tt 4.2480e+05} & {\tt 7.6889e+00} &  {\tt 1.0431e+14} \\
$60$   & $30$ & {\tt 5.0e-11} & {\tt 1.7761e-02} & {\tt 7.2223e+04} & {\tt 2.4873e+00} &  {\tt 1.3124e+11} \\
$120$  & $30$ & {\tt 5.0e-12} & {\tt 7.2560e-02} & {\tt 3.7582e+03} & {\tt 4.8933e-01} &  {\tt 7.4516e+05} \\
\end{tabular}
\caption{{\sc Information for $\ell=64$ compressed kernels.}
Note that the min/max values for $|a_{64,j}(r_1,240)|$ have been
computed with \eqref{eq:afromscaled}. Values computed with
the asymptotic expansion \eqref{eq:mainresult} are the same to 
about 3 digits of relative accuracy.
\label{tab:ell64}
}
\end{table*}

\begin{figure}[htp]
\begin{center}
\subfigure[Pan out.]{
\includegraphics[width=7.5cm]{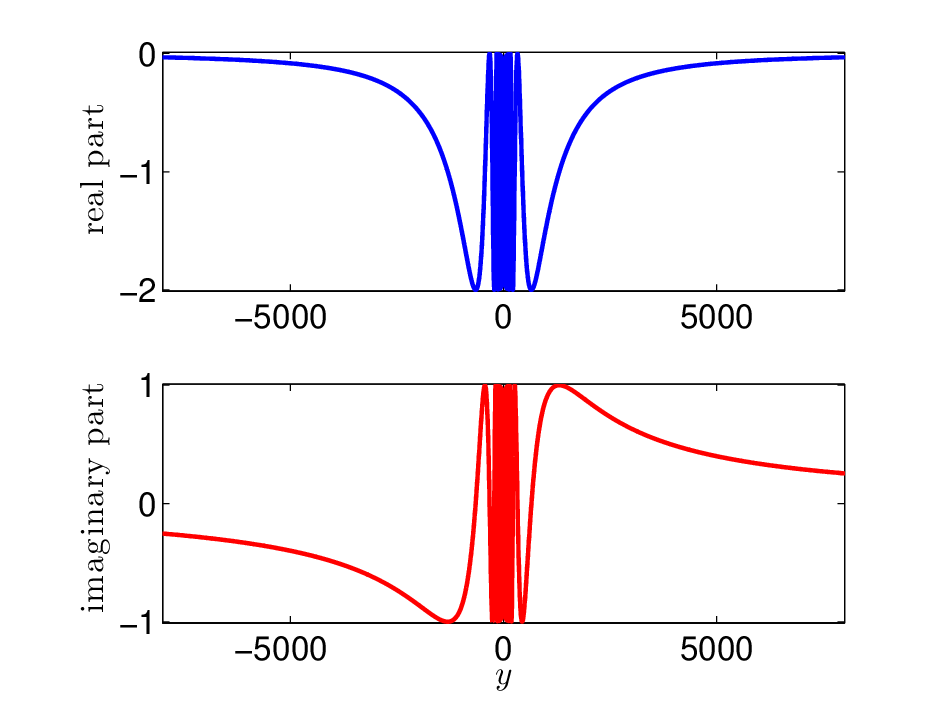}
}
\subfigure[Zoom in.]{
\includegraphics[width=7.5cm]{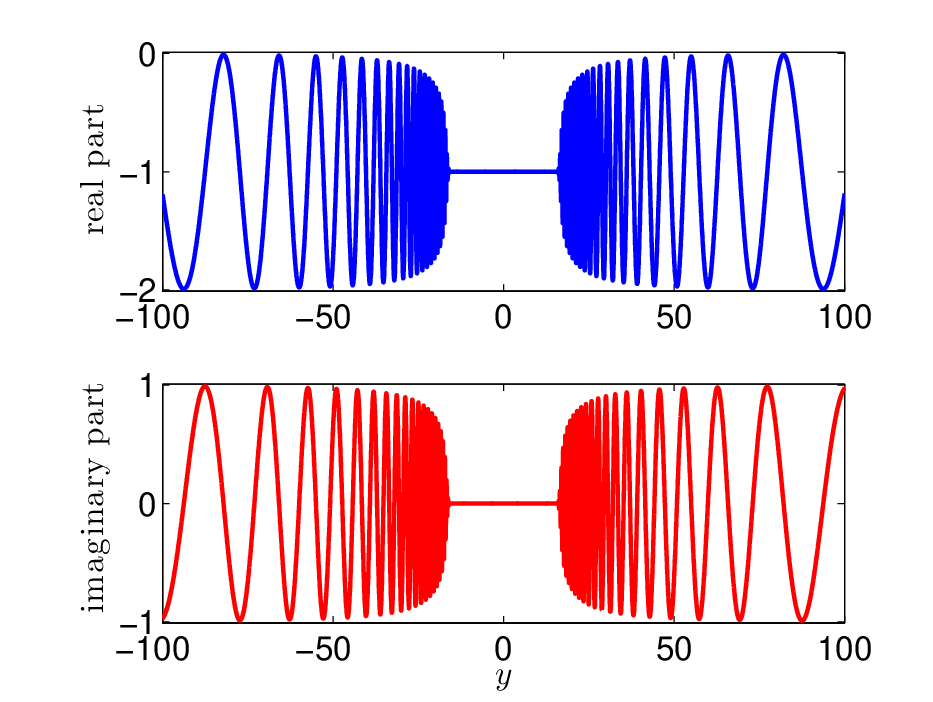}
}
\end{center}
\caption{Teleportation kernel $\widehat{\Phi}_{256}(\mathrm{i}y,15,240)$.
The upper plots depict Re$\widehat{\Phi}_{256}(\mathrm{i}y,15,240)$
and the bottom plots Im$\widehat{\Phi}_{256}(\mathrm{i}y,15,240)$.}
\label{fig:ell256kernel}
\end{figure}

As a more extreme example, consider the teleportation
kernel $\widehat{\Phi}_{256}(\mathrm{i}y,15,240)$ 
shown in Fig.~\ref{fig:ell256kernel}. Based on \eqref{eq:mainresult}, 
we expect the residues $a_{256,j}(15,240)$ to vary in modulus over 
64 orders of magnitude, with the largest residue $\gtrsim 10^{67}$ 
in magnitude. Therefore, a table based on storage of the exact 
poles and residues in double precision format likely yields no
digits of accuracy. Nevertheless, we have found that a compressed
kernel with $d=72$ poles, when stored in double precision format,
corresponds to $\varepsilon \lesssim$ {\tt 5.0e-08}; moreover, storage
of a $d=79$ compressed kernel in quadrupole precision format
corresponds to $\varepsilon \lesssim$ {\tt 5.0e-11}.

\subsection{Empirical scalings} \label{sec:num_scalings}

We have carried out 1183 independent kernel compressions over the
parameters $\left(\ell, d, \varepsilon, r_2/r_1\right)$, and with
these explore relevant scalings for our approximations.
Figures~\ref{fig:err_vs_pole} to \ref{fig:r_vs_res} summarize the
results. Unless stated otherwise, all experiments in this 
subsection vary $r_1$ with $r_2 = 240$ fixed. Results 
for $r_1 = 15$ fixed and varying $r_2$ are in most cases 
qualitatively similar. We are primarily interested in the first 
approach to varying $r_2/r_1$, since it comports with 
teleportation as a technique for asymptotic waveform evaluation
\cite{BFL}, although typically $r_2 \gg 240$ in that context.
Each compression has been computed in quadruple precision arithmetic,
with the resulting pole locations and residues stored in both 
double and quadruple precision formats. For smaller 
$\ell$ values we have used the double precision format without loss 
of accuracy. However, in some cases (e.g.~$\ell = 128,256$) below,
we have found it necessary to use the quadruple precision format
to retain high accuracy.

\subsubsection{Scaling of
compressed-kernel pole number $d$ with compression 
error $\varepsilon$}  
\label{subsubsec:error_scaling}

\begin{figure}[htp]
\begin{center}
\includegraphics[width=7.55cm]{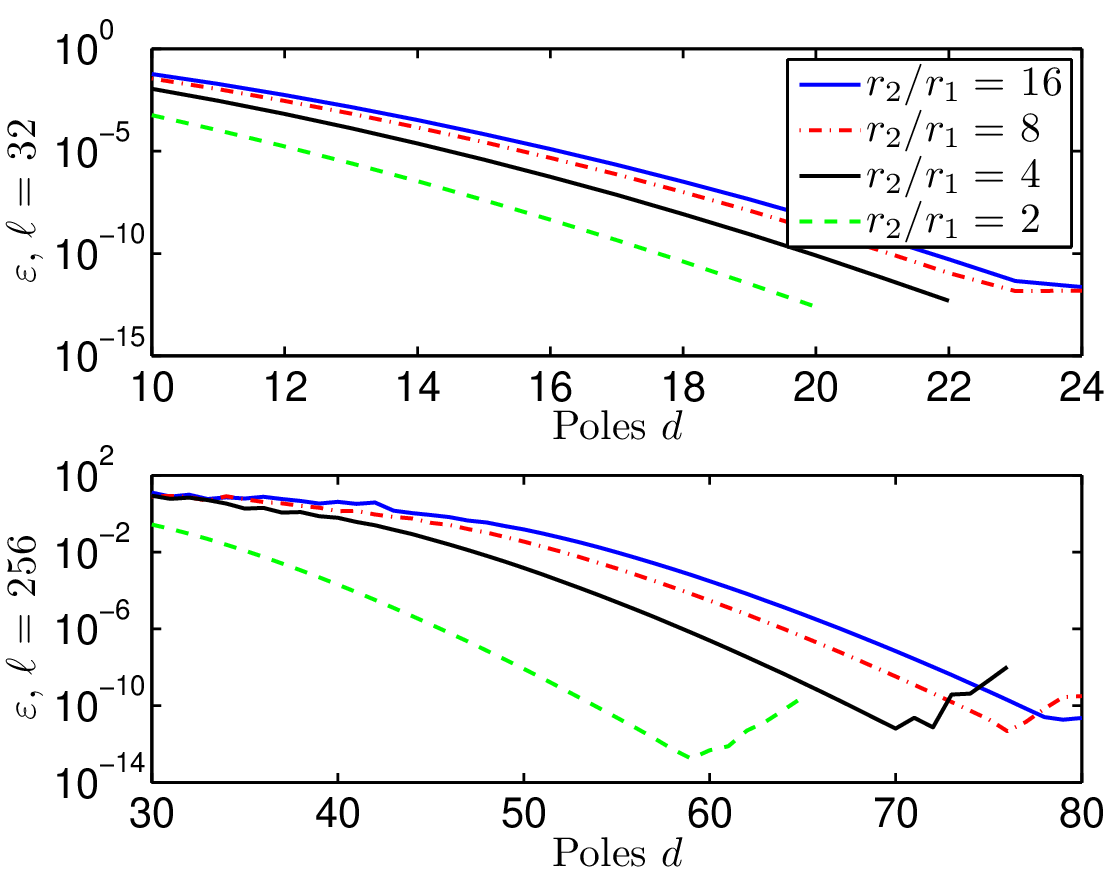}
\includegraphics[width=7cm]{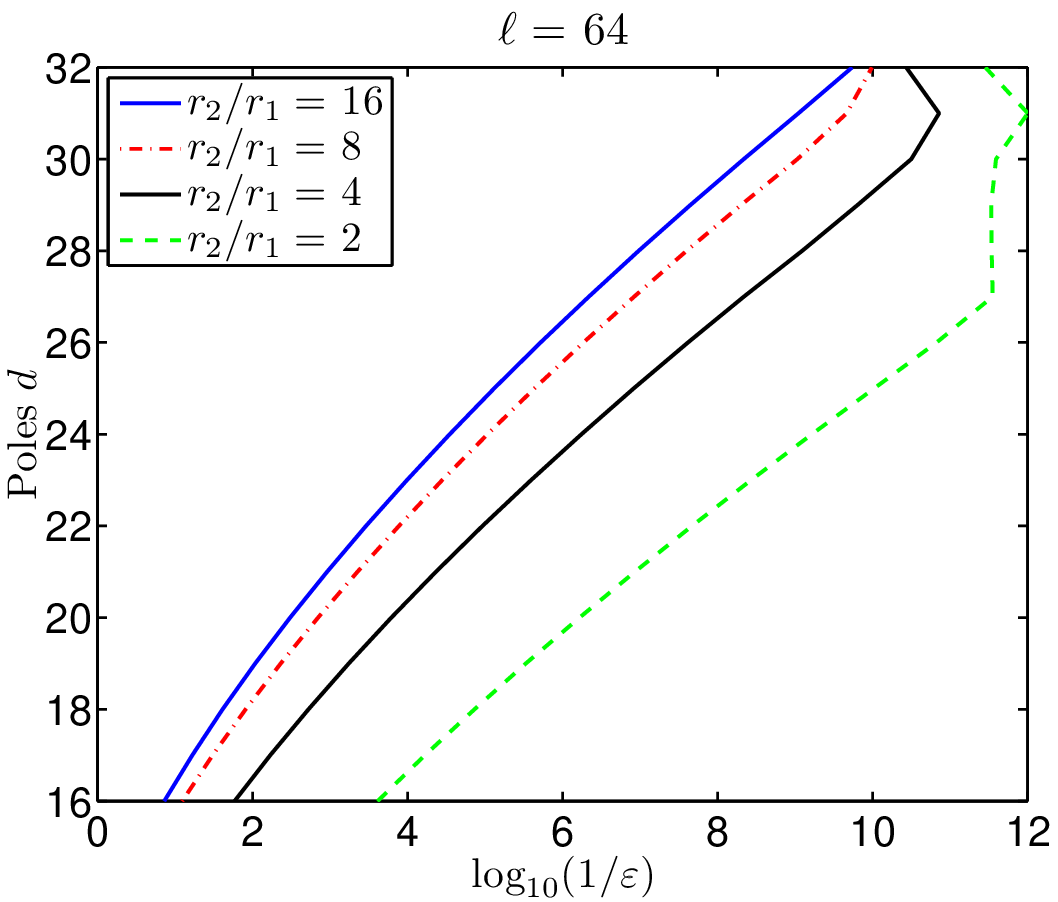}
\end{center}
\caption{Exponential decay of compression error 
$\varepsilon$ with pole number $d$. Plots for
other $\ell$ values are qualitatively similar.}
\label{fig:err_vs_pole}
\end{figure}

Figure~\ref{fig:err_vs_pole}(left) confirms the exponential decay 
of the approximation suggested by Lemma~\ref{lem:pole_approx}. 
These plots indicate that for fixed values of $\ell$ and $r_2/r_1$, 
the number $d$ of approximating 
poles scales like
\begin{align}\label{eq:scaling_d_with_eps}
d = O\big(\log(1/\varepsilon)\big).
\end{align}

\subsubsection{Scaling of compressed-kernel pole
number $d$ with spherical-harmonic index $\ell$}

\begin{figure}[htp]
\begin{center}
\includegraphics[width=6.5cm]{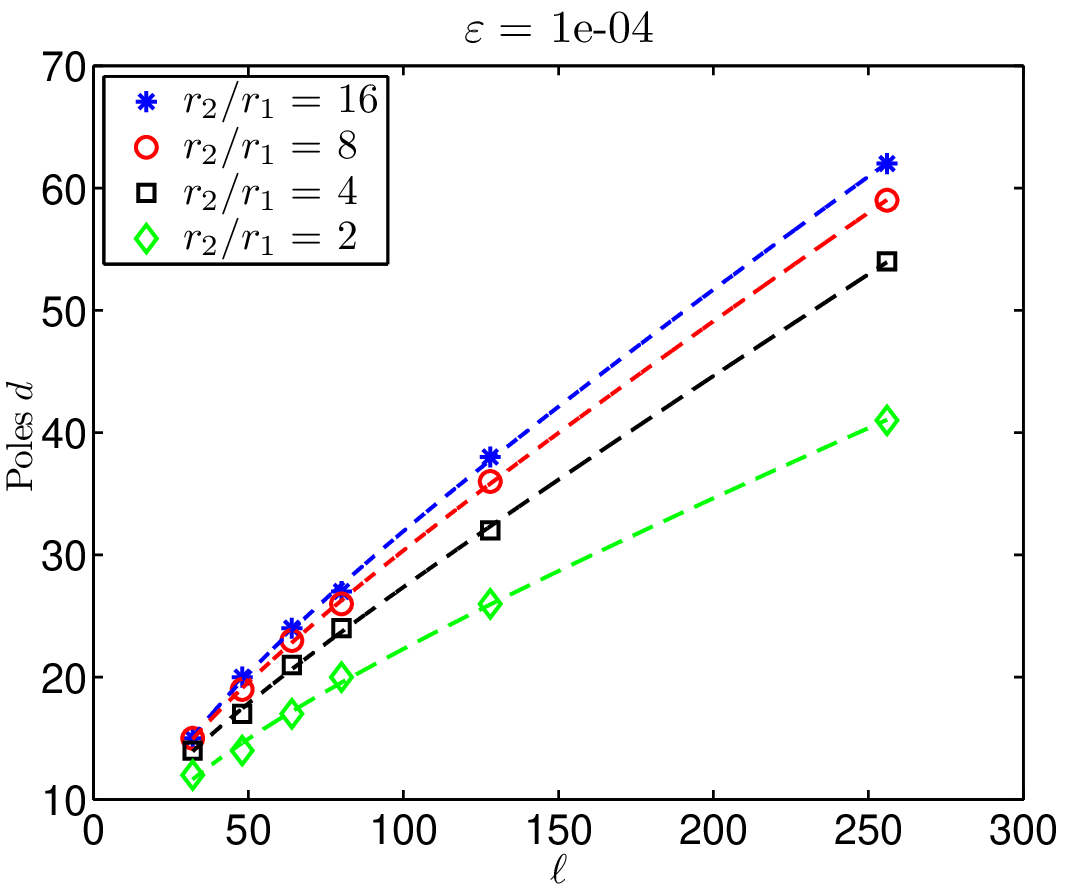}
\includegraphics[width=6.5cm]{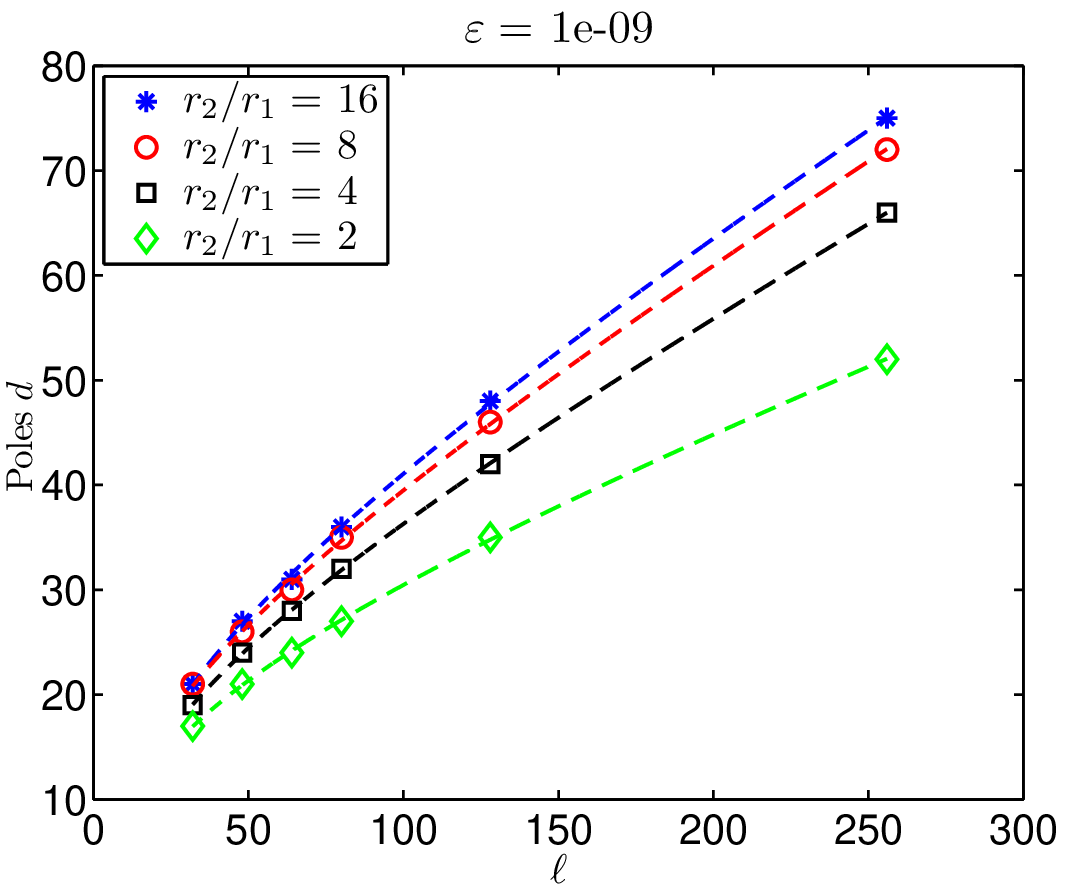}
\end{center}
\caption{Growth of pole number $d$ with 
spherical-harmonic index $\ell$. 
Plots for other values of $\varepsilon$ 
are qualitatively similar. Fits described by
\eqref{eq:scaling_d_with_ell} correspond to the 
dashed lines.}
\label{fig:ell_vs_poles}
\end{figure}

Figure~\ref{fig:ell_vs_poles} depicts the growth of the 
compressed-kernel pole number $d$ with spherical-harmonic 
index $\ell$. For fixed values of $r_1$ and $r_2$, 
the number of approximating poles is best described by the fit
\begin{align}\label{eq:scaling_d_with_ell}
d = a_1 + a_2 \log\left(\ell\right) + a_3 \ell.
\end{align}
As judged by standard goodness-of-fit statistics, 
Eq.~\eqref{eq:scaling_d_with_ell} models the growth
of $d$ better than solely $\log$, linear, or $\log^2$ scaling.
Formally then, $d$ would appear to be $O(\ell)$
for large $\ell$. However, for 
small and medium $\ell$ the logarithmic 
term dominates, and, furthermore, the coefficient $a_2$ is typically 
small. For example, the fits shown in 
Fig.~\ref{fig:ell_vs_poles}(left) are as follows.
\begin{verbatim}
    r2/r1 = 16: a1 = -7.814, a2 = 5.112, a3 = 0.1621
    r2/r1 =  8: a1 = -4.096, a2 = 4.017, a3 = 0.1596 
    r2/r1 =  4: a1 = 0.1785, a2 = 2.537, a3 = 0.1551 
    r2/r1 =  2: a1 = -3.449, a2 = 3.435, a3 = 0.0994 
\end{verbatim}
The scaling~\eqref{eq:scaling_d_with_ell} is at odds
with the conjecture made in Ref~\cite{BFL} just after Eq.~(40)
of that reference.

\subsubsection{Scaling of largest residue $\max_j |\gamma_{\ell j} |$ with 
compression error $\varepsilon$}
\begin{figure}[htp]
\begin{center}
\includegraphics[width=6.5cm]{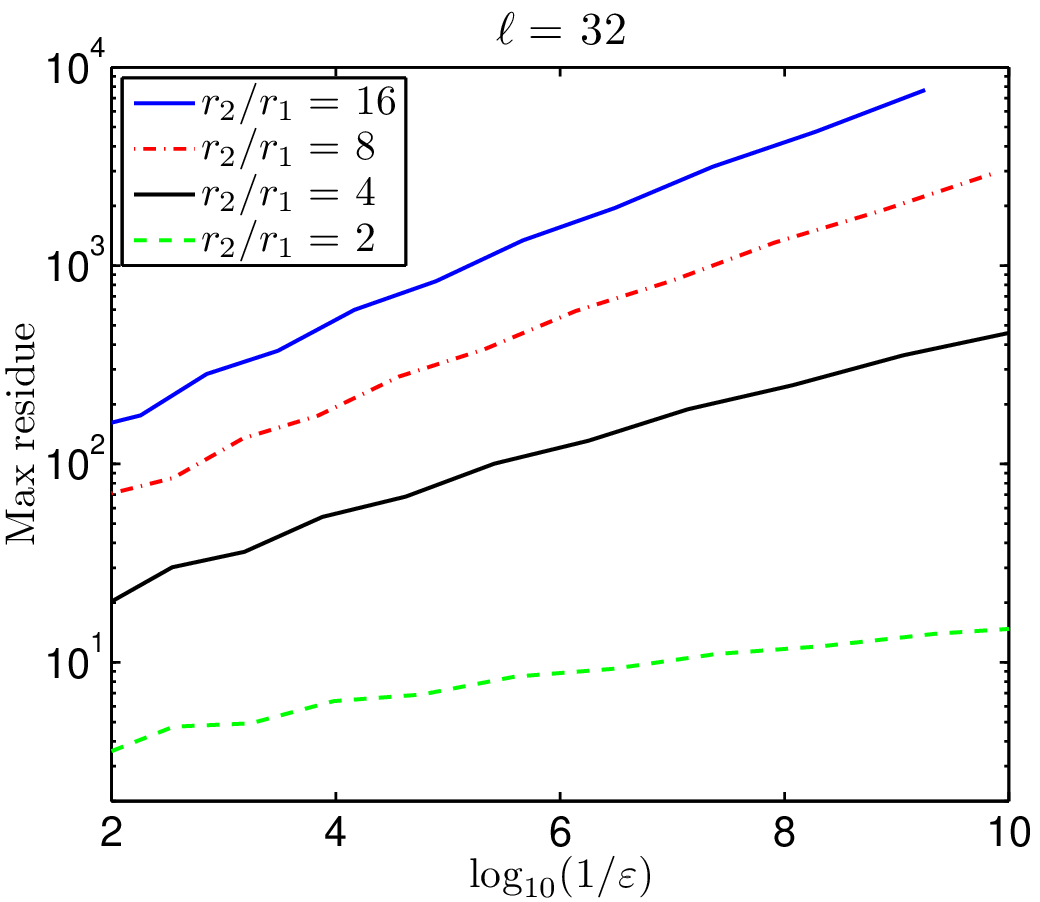}
\includegraphics[width=6.5cm]{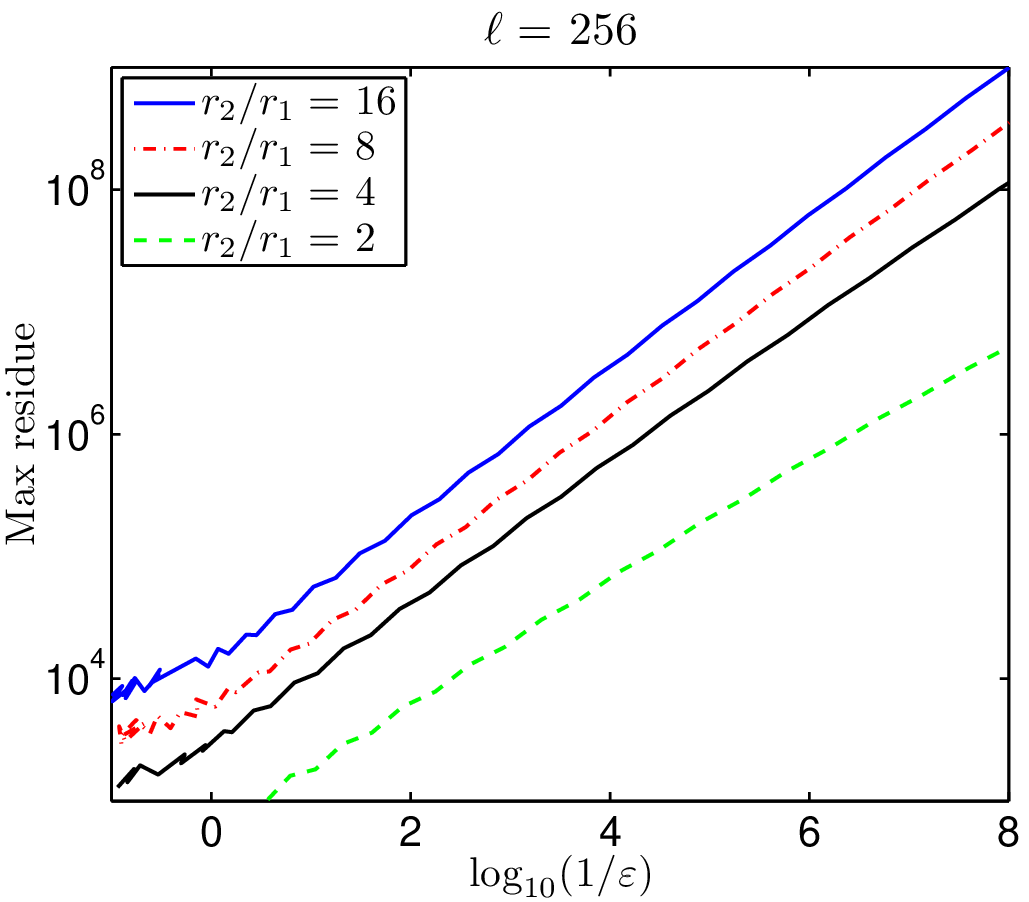}
\end{center}
\caption{Growth of compressed-kernel
largest residue $\max_j |\gamma_{\ell j}|$ 
with compression error $\varepsilon$.}
\label{fig:pole_vs_resmax}
\end{figure}

Figure~\ref{fig:pole_vs_resmax} depicts the growth of the 
compressed-kernel maximum residue (in modulus) with
the approximation. These plots indicate that for
fixed values of $\ell$, $r_1$, and $r_2$ 
the largest residue scales like
\begin{align}
\max_j |\gamma_{\ell j} | = O\big(\log(1/\varepsilon)\big).
\end{align}
We anticipate that $\max_j |\gamma_{\ell j} |$ approaches
$\max_j | a_{\ell j} |$ as $d \rightarrow \ell$, since
the best $\ell$-pole approximation of the kernel would be
the kernel itself. Due to the high precision required we are 
unable to probe the regime $d \approx \ell$ in the large-$\ell$ 
limit.

\subsubsection{Scaling of largest residue $\max_j |\gamma_{\ell j} |$ with 
spherical-harmonic index $\ell$}

\begin{figure}[htp]
\begin{center}
\includegraphics[width=6.5cm]{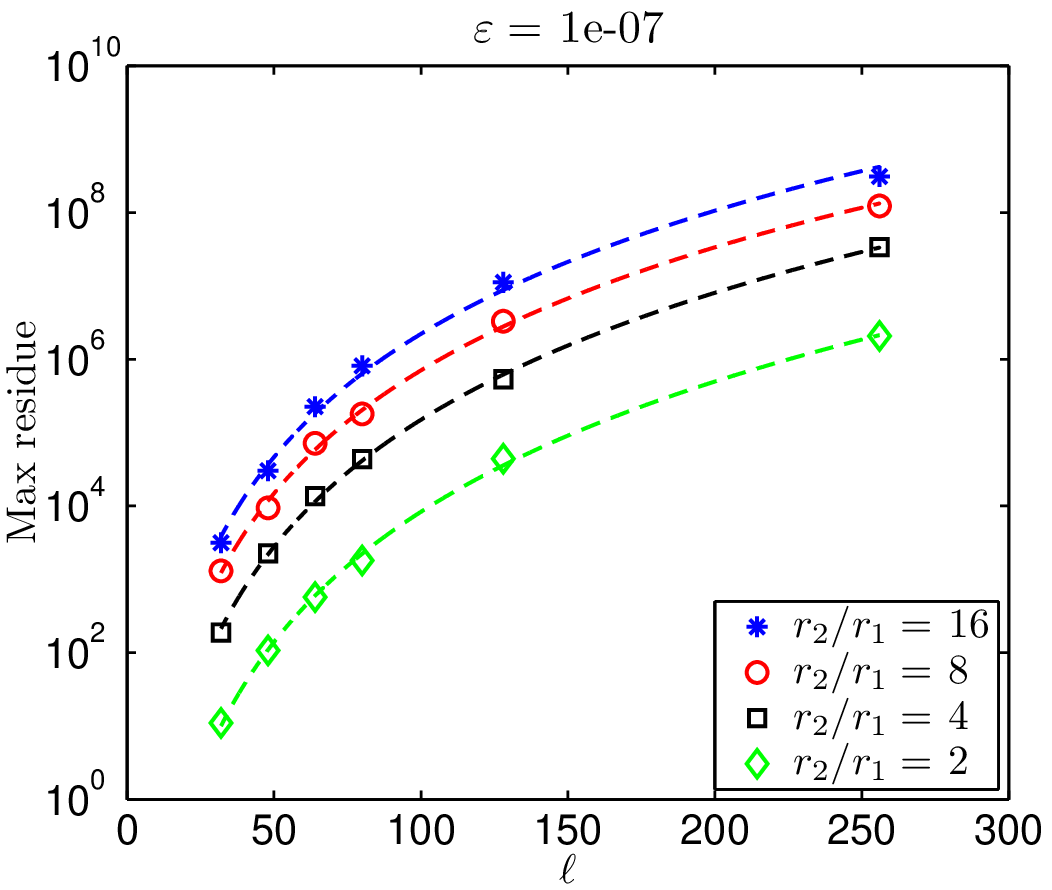}
\includegraphics[width=7.6cm]{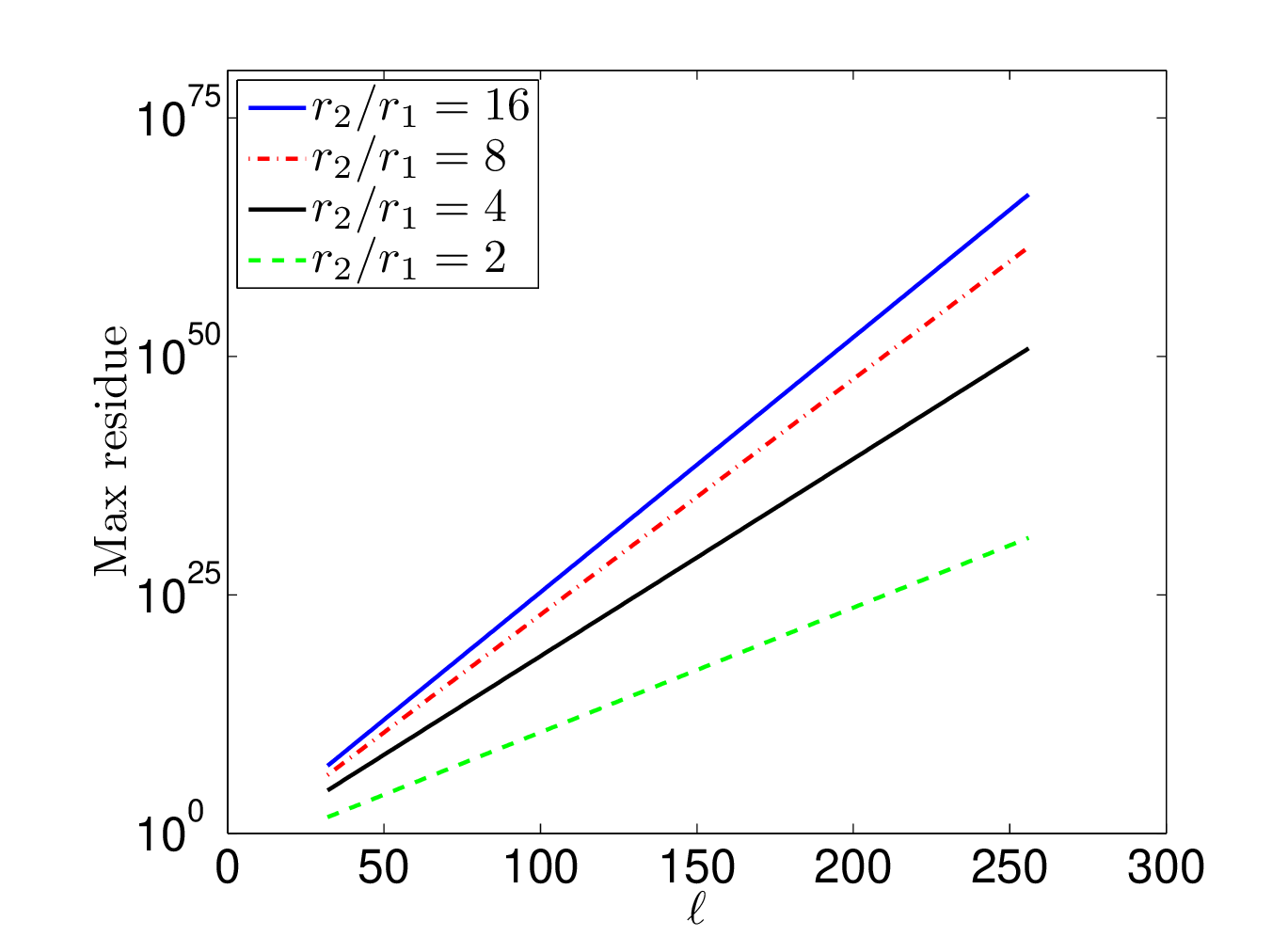}
\end{center}
\caption{Growth of compressed-kernel largest 
residue $\max_j |\gamma_{\ell j}|$
with spherical-harmonic index $\ell$. 
Results for four families of $\varepsilon = 10^{-7}$ 
kernels are shown(left), with the dashed 
lines showing fits to the data as described by
Eq.~\eqref{eq:maxgam_from_ell}.
For comparison, the corresponding results for the
exact-kernel $\max_j |a_{\ell j}|$ 
are also shown(right).
}
\label{fig:ell_vs_res}
\end{figure}

Figure~\ref{fig:ell_vs_res}(left) depicts the growth of 
the compressed-kernel maximum residue (in modulus) 
with spherical-harmonic index 
$\ell$. For fixed values of $r_1$, $r_2$, and $\varepsilon$ 
the data is well modeled by
\begin{align}\label{eq:maxgam_from_ell}
\max_j |\gamma_{\ell j} | = 
\exp\big( a_1 \log (\ell) + a_2 \big).
\end{align}
The fits shown in Fig.~\ref{fig:ell_vs_res}(left) are as 
follows.
\begin{verbatim}
    r2/r1 = 16: a1 = 5.58, a2 = -11.09 
    r2/r1 =  8: a1 = 5.57, a2 = -12.20
    r2/r1 =  4: a1 = 5.75, a2 = -14.57
    r2/r1 =  2: a1 = 5.80, a2 = -18.14 
\end{verbatim}
Figure~\ref{fig:ell_vs_res}(right) depicts the
growth of the exact-kernel
maximum residue (in modulus), as computed by \eqref{eq:mainresult}.

\subsubsection{Scaling of largest residue $\max_j |\gamma_{\ell j} |$ with 
$r_2/r_1$}

\begin{figure}[htp]
\begin{center}
\subfigure{
\includegraphics[width=7cm]{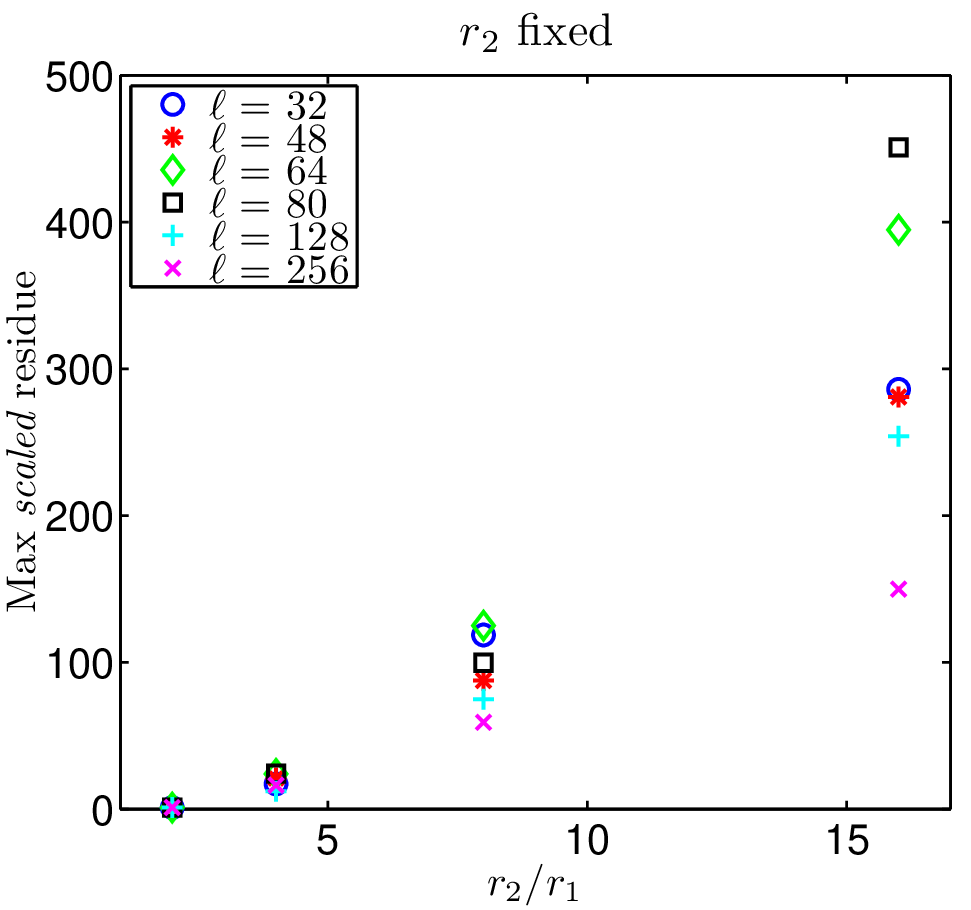}
\includegraphics[width=7cm]{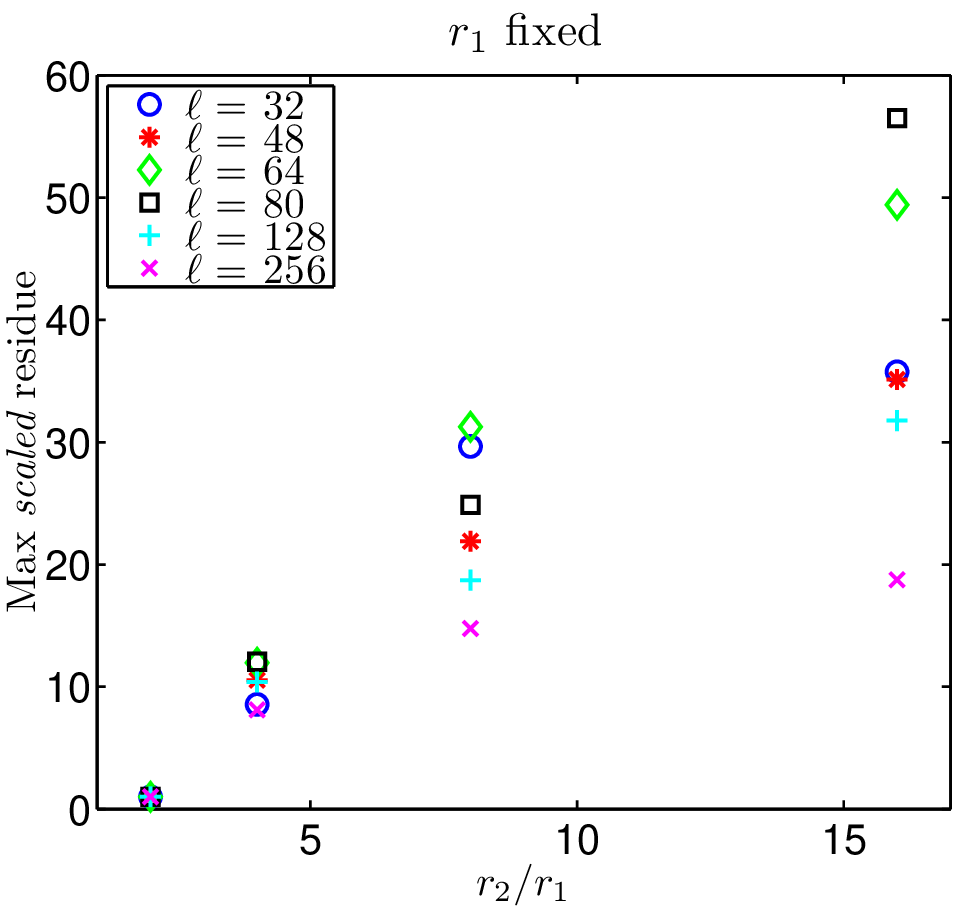}
}
\end{center}
\caption{Growth of compressed kernel's largest residue with 
pole number $r_2/r_1$. 
For plotting convenience a scaled maximum residue is shown. The scaling factor is 
chosen such that the maximum residue is $1$ at $r_2/r_1 = 2$, thereby ensuring
similar residue magnitudes among the different $\ell$ kernels.}
\label{fig:r_vs_res}
\end{figure}

Figure~\ref{fig:r_vs_res} depicts the growth of the 
compressed-kernel maximum 
residue (in modulus) with $r_2/r_1$. We consider two cases: 
(i) $r_2$ fixed (left) and (ii) $r_1$ fixed (right). Notice 
that the maximum residue grows more quickly as 
$r_1 \rightarrow 0$ than it does as $r_2 \rightarrow \infty$.
This observation is expected. Indeed, for case (i) the effective 
potential $\ell(\ell+1)r_1^{-2}$ at $r_1$ grows without bound as 
$r_1 \rightarrow 0$. Whence increasing $r_2/r_1$ corresponds 
to propagation from a region of increasingly large potential.
However, for case (ii) similar considerations show that 
increasing $r_2/r_1$ corresponds to propagation into a region 
of increasingly small potential.
Regardless, the limiting case 
$r_2 \rightarrow \infty$ does not appear problematic for 
those cases considered. Furthermore, 
Ref.~\cite{BFL} has achieved high accuracy kernel compressions 
for $r_2 \simeq 10^{15}$ (albeit for low-$\ell$ and
the wave equations describing gravitational perturbations).

\subsection{Pulse teleportation}\label{subsec:pulse}
This subsection presents an experiment similar 
to the one carried out in \cite{BFL} for the Regge-Wheeler 
equation (see the conclusion).
We take $\ell = 64$ with Gaussian initial data:
\begin{equation}
\Psi(0,r) = e^{-(r-8)^2},\qquad 
\partial_t\Psi(0,r) = -\partial_r\Psi(0,r),
\end{equation}
where we have suppressed the $(\ell=64,m)$ indices on $\Psi$. 
Using a multidomain nodal Chebyshev method (described in \cite{FHL}), 
we perform five separate evolutions on domains with outer boundaries 
$b=15,30,60,120,240$. We have respectively used 8, 16, 30, 60, and 
120 subintervals of uniform size, and in each case with 42 
Chebyshev-Lobatto points per subinterval. Therefore, the spatial
resolution for each evolution is comparable to the others.
Evolutions are performed by the classical 4-stage explicit Runge 
Kutta method with timestep $\Delta t \simeq$ 
{\tt 4.0461e-05}. For each evolution the inner boundary 
is $a = 2$, at which we have enforced a Sommerfeld boundary 
condition $\partial_t\Psi-\partial_r\Psi = 0$. For all
choices of outer boundary $b$ we adopt the Laplace convolution
radiation boundary conditions (RBC) 
based on (compressed kernels for) the time-domain kernel 
$\Omega_{64}(t,b)$ (see \cite{BFL} for more detail). Tables for 
the $b = 15,20,60,120,240$ RBC respectively have 19, 19, 19, 18, 
and 17 poles, with each table computed in quadruple precision to 
satisfy the tolerance $\varepsilon = 10^{-15}$. 

\begin{figure}[htp]
\begin{center}
\includegraphics[width=11cm]{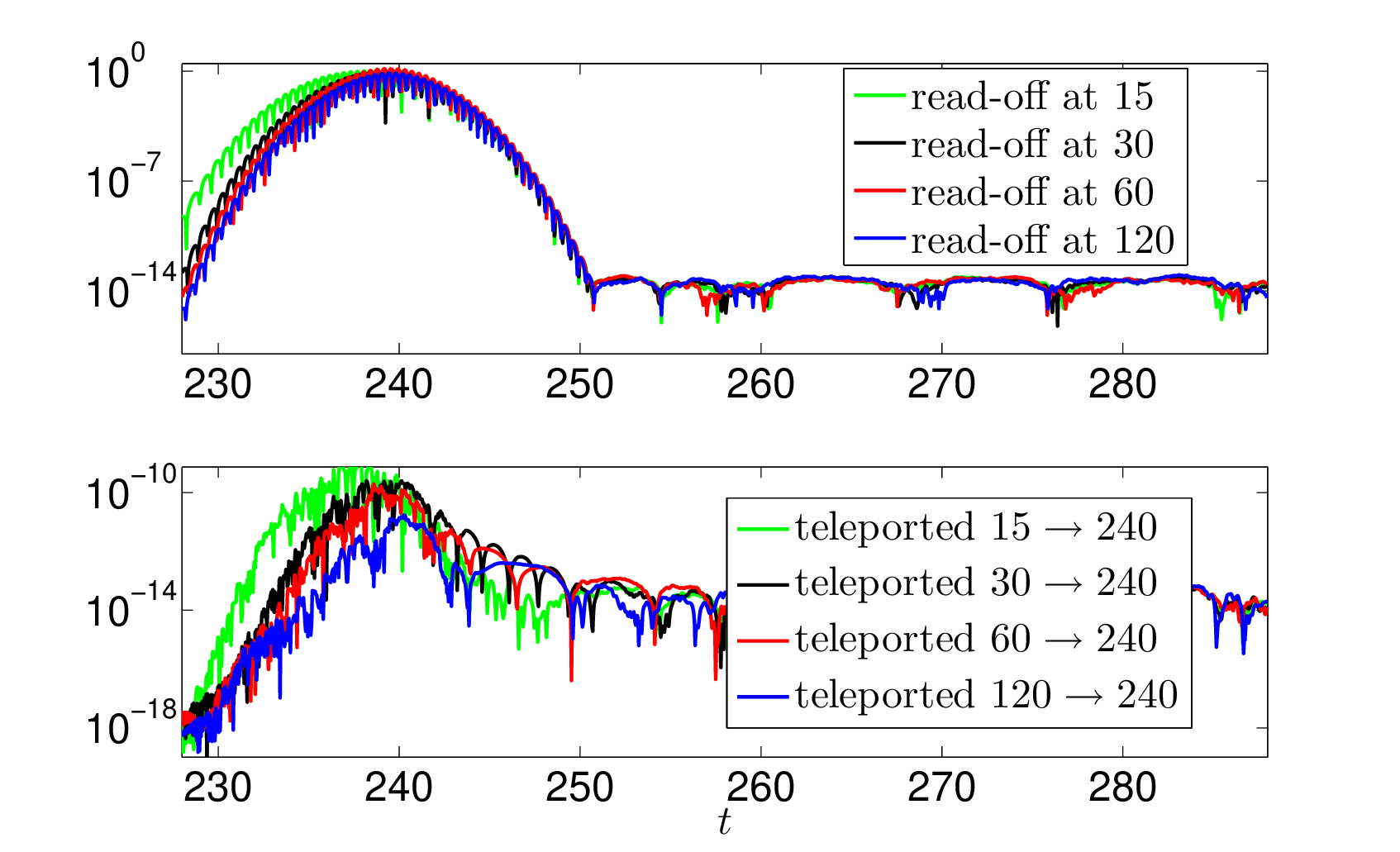}
\end{center}
\caption{Teleportation of an $\ell=64$ Gaussian pulse.
Here all errors have been computed against the reference
signal read-off at $r = 240$. The top panel then compares
the reference signal with the signals read-off
at $r=15,30,60,120$. The bottom panel compares the reference
signal with the signals obtain through teleportation 
from $r_1=15,30,60,120$ to $r_2=240$.}
\label{fig:pulse64}
\end{figure}

In all five simulations the field $\Psi(t,b)$ is recorded as a time 
series at the boundary $b$, and in all cases but $b=240$
we "teleport" the field from $r_1 = b$ to $r_2 = 240$.
The compressed teleportation kernels
$\Xi_{64}(t,b,240) \simeq \Phi_{64}(t,b,240)$ have been described in
the previous subsection. For the $b=240$
simulation we simply record the field at the boundary, with
this record then serving as a reference time series.
We account for time delays by starting all recorded times
series (whether read off or teleported) at time $b - 12$.
The top panel in Figure \ref{fig:pulse64} plots 
the errors in the waveforms recorded at the different $b$
boundaries as compared to the reference $b=240$ waveform;
as expected the systematic errors are large.
The bottom panel plots the errors in the
$r_1 = b \rightarrow r_2 = 240$ teleported
time series relative to the reference time series.

\subsection{Asymptotic-waveform evaluation}\label{subsec:awe}
We now consider asymptotic-waveform evaluation (AWE) by
teleportation of a pure $\ell=2$
multipole signal recorded from $r_1 = 10$ to 
$r_2 = \infty$. From \eqref{eq:fdkernel_poles} the 
$\ell = 2$ frequency-domain teleportation kernel is
\begin{align}
\widehat{\Phi}_2(s,r_1,r_2) = 
\frac{a_{2 1}(r_1,r_2)}{s - z_{-}/r_1} 
+ \frac{a_{2 2}(r_1,r_2)}{s - z_{+}/r_1} , \quad 
z_{\pm} = -\frac{3}{2} \pm \mathrm{i}\frac{\sqrt{3}}{2},
\end{align}
where $W_2(z_{\pm}) = 0$ 
and $b_{21} = z_-$, $b_{22}=z_+$. 
The AWE residues can be found through direct evaluation 
of \eqref{eq:residues}, giving
\begin{align} \label{eq:ell2awe}
\Phi_2(t,10,\infty) = 
\tilde{a} \exp\left(\frac{z_{-}}{10}t \right) + 
\tilde{a}^* \exp\left(\frac{z_{+}}{10}t \right), \quad 
\tilde{a} = \mathrm{i} \frac{\sqrt{3}}{30} z_{-}^2, 
\end{align}
as the relevant time-domain AWE kernel and where we 
make use\footnote{For even $\ell$ the residues 
come in conjugate pairs, as seen
from \eqref{eq:residues} and 
the fact that the roots $b_{\ell j}$ also come in 
conjugate pairs.}
of $a_{2 2} = a_{2 1}^*$.

We now use the AWE kernel \eqref{eq:ell2awe} to recover 
the asymptotic solution from a recorded time-series.
To directly compute errors, we shall consider an exact, 
closed-form solution. From Eq.~\eqref{eq:ellmult_fd} 
the general $\ell=2$ outgoing solution
to Eq.~\eqref{eq:radialwave} is
\begin{align} \label{eq:outgoingEll2}
\Psi(t,r) = f''(t-r) 
                + \frac{3}{r} f'(t-r) 
                + \frac{3}{r^2} f(t-r),
\end{align}
where $f(u)$ is an underlying function of retarded 
time $u=t-r$, the prime indicating 
differentiation in argument, and we have
suppressed the $(\ell=2,m)$ indices on $\Psi$.
The specification 
\begin{align}\label{eq:sinegauss}
f(u) = \sin \left[f_0 \left(u-u_0\right) \right] 
\mathrm{e}^{-c\left(u-u_0\right)^2},
\end{align}
determines a purely outgoing multipole solution whose 
asymptotic signal is
\begin{align} \label{eq:AWEsig_ell2}
\Psi_{\infty}(T) \equiv 
\left[ -4 c f_0 T \cos\left(f_0 T \right) +
\left(-2 c - f_0^2+4 c^2 T^2\right) 
\sin \left(f_0 T\right)  \right]
\mathrm{e}^{-cT^2}, \quad T \equiv t - r_1 - u_0.
\end{align}
Here $c$ characterizes the solution's spatial extent,
$f_0$ its "central" frequency, and $u_0$ its offset.
To obtain \eqref{eq:AWEsig_ell2}, we have 
adjusted for the infinite time delay for the signal to reach 
$r_2=\infty$. Indeed, the signal at $(t_1,r_1)$ reaches 
$r_2$ at time $t_1+(r_2-r_1)$; both correspond to the retarded time
$u_1=t_1-r_1$ (the combination in $T$) even for $r_2=\infty$.
For all experiments we choose 
$f_0 = 1$, $c = 2.5$, and $u_0 = -6$. To distinguish 
from exact solutions, we append a subscript $h$ to any quantity 
obtained with our multidomain nodal Chebyshev solver.

\begin{figure}[htp]
\centering
\subfigure[Numerically generated time-series]{
  \includegraphics[width=7.3cm]{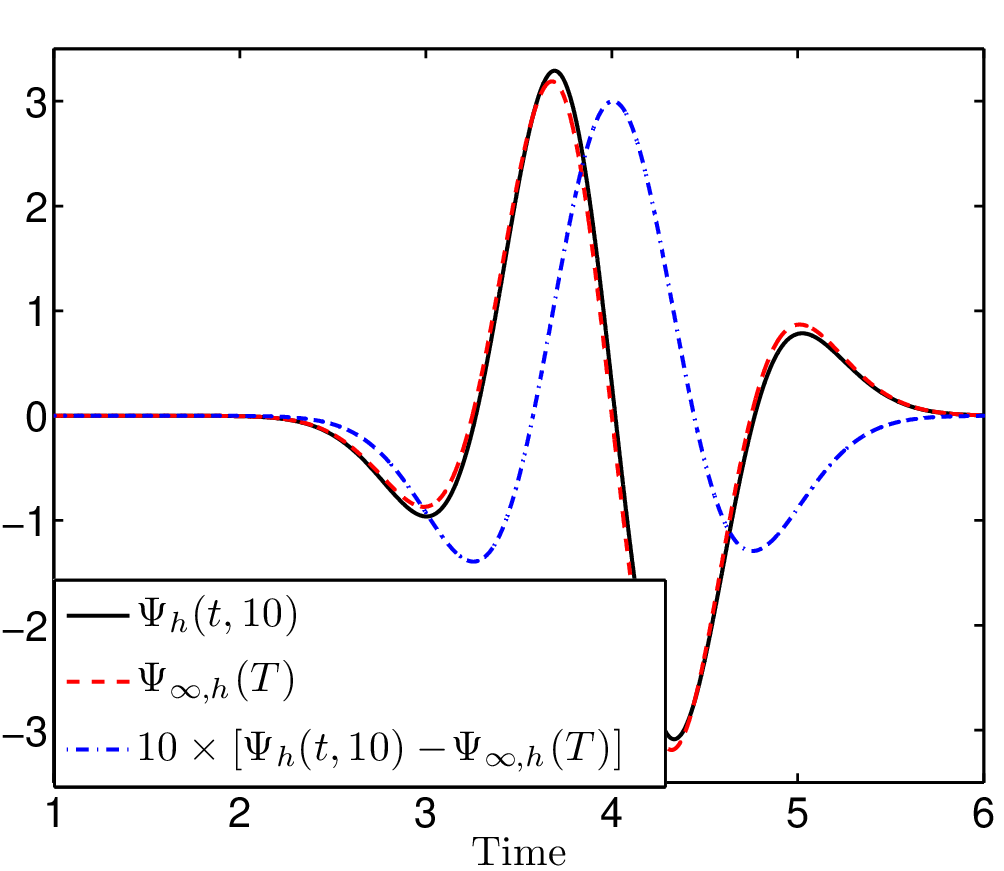}
  \label{fig:psiAwe}}
\subfigure[Time-series error]{
  \includegraphics[width=7.3cm]{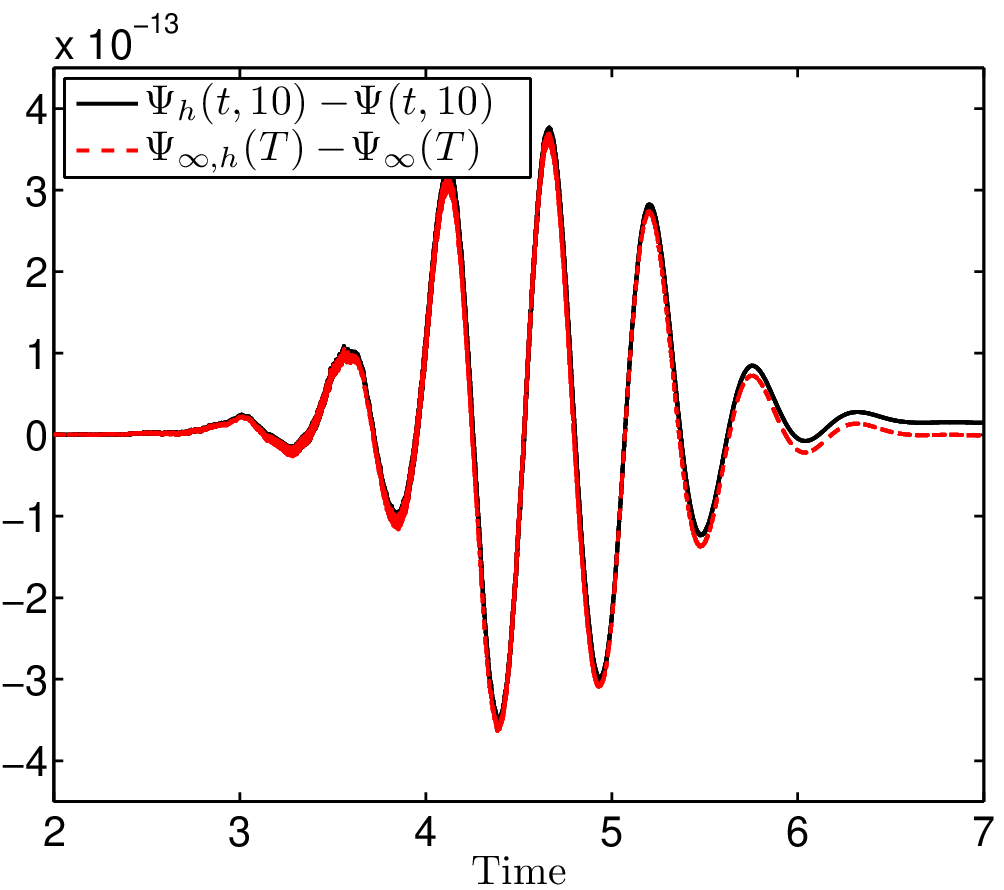}
  \label{fig:psiAwe_err}}
\caption{The left plot depicts the $r_1 = 10$ recorded data 
(solid black line) and its teleportation to $r_2 = \infty$ 
(dashed red line). The difference between these 
time-series (dash-dot blue line) characterizes the 
systematic error engendered by approximation of 
the true asymptotic signal by the  finite-radial value one.
For visual assistance, we have multiplied this error by $10$. 
The right figure depicts the numerical errors which have
been computed relative to the exact time-series.}
\label{fig:psiAwe_exp}
\end{figure}

The numerical setup is similar that of Sec.~\ref{subsec:pulse},
except that we now choose 22 Chebyshev-Lobatto points on 8 subintervals
and $\Delta t \simeq$ {\tt 5.3949e-05}. Initial data 
is found from Eq.~\eqref{eq:outgoingEll2} and its spatial 
and temporal derivatives, all evaluated at $t=0$. The 
signal $\Psi_h(t,10)$ recorded at $r_1 = 10$ is the 
solid black line in Fig.~\ref{fig:psiAwe}. 
With this data we compute $\Psi_{\infty,h}\left(T\right)$,
the dashed red line in Fig.~\ref{fig:psiAwe}, from
the convolution \eqref{eq:td_teleportation} with kernel 
\eqref{eq:ell2awe}. In numerical studies, often the 
outermost recorded signal is taken to be the asymptotic 
one. The dash-dot blue line in Fig.~\ref{fig:psiAwe} plots
the {\em systematic} error
$\Psi_{\infty,h}\left(T\right) -\Psi_h(t,10)$,
with multiplication by 10 as a visual aid.
Notice that the systematic error is $O\left(1/r_1\right)$, 
in accord with Eq.~\eqref{eq:outgoingEll2}. 

The difference between the 
exact asymptotic signal $\Psi_{\infty}(T)$ and our numerically 
teleported one is the dashed red line in Fig.~\ref{fig:psiAwe_err}. 
As the exact AWE kernel has been used with
teleportation convolution that is accurately evaluated for our 
choice of $\Delta t$, these errors 
must stem from the numerical solver.
This expectation is confirmed by the solid black line in 
Fig.~\ref{fig:psiAwe_err} which plots the error 
$\Psi(t,10) - \Psi_h(t,10)$.

\begin{figure}[htp]
\begin{center}
\includegraphics[width=8cm]{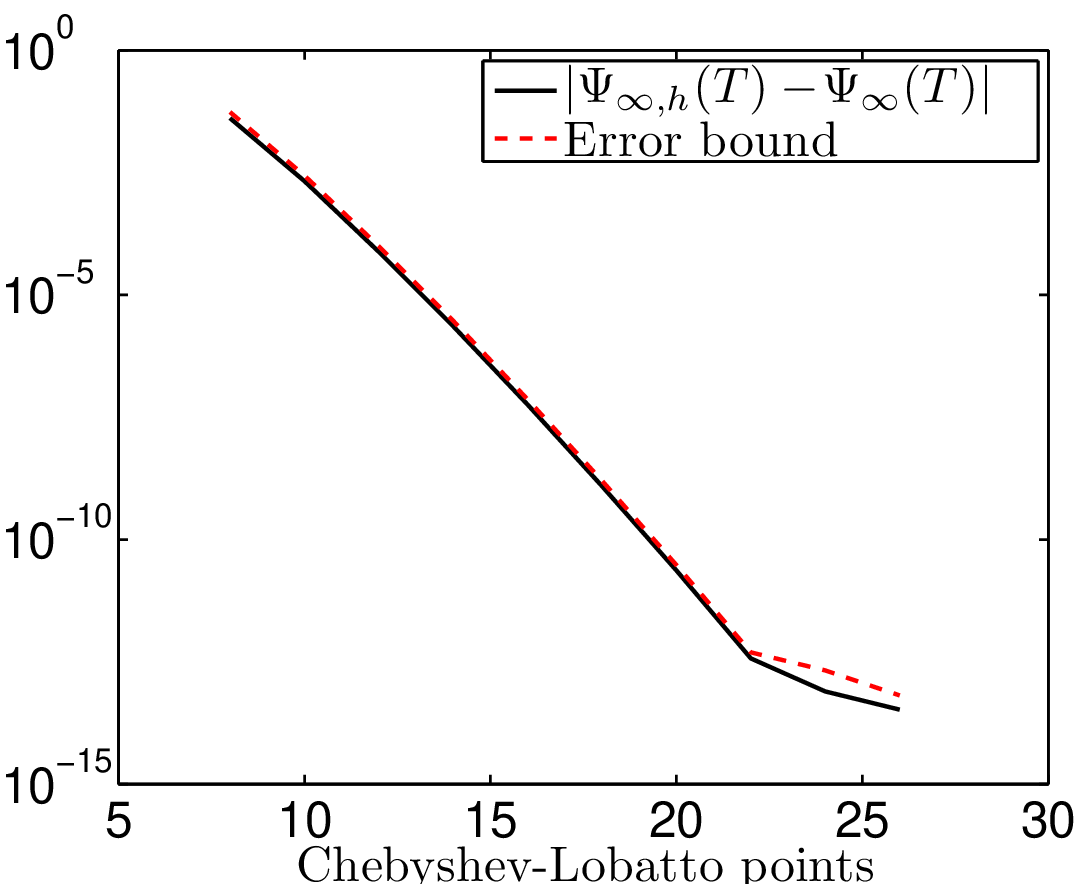}
\end{center}
\caption{The exact (solid black line) and estimated (dashed red line) 
error from Eq.~\eqref{eq:errbound}. To form this bound we 
compute $\big\| \Phi_2(\cdot,10,\infty) \big\|_{L_{\infty}} = 
\frac{3}{10}$ from Eq.~\eqref{eq:kernLinf}.
}
\label{fig:psiAwe_bound}
\end{figure}

When teleportation is carried out to a large radial value, 
performance of a convergence test on a commensurately larger 
computational domain becomes unfeasible. Such is the case
here, but Eq.~\eqref{eq:errbound} from the appendix
provides a useful error bound for 
$\Psi_{\infty}\left(T\right) -\Psi_{\infty,h}(T)$ in 
terms of $\Psi\left(t,10\right) -\Psi_h(t,10)$, where the
latter error can be estimated through a convergence study.
The dashed red line in Fig.~\ref{fig:psiAwe_bound}
plots this error bound versus the number of Chebyshev-Lobatto 
points used in our computation and sharply bounds the 
solid black line which depicts the exact error. 

\section{Concluding remarks}

For solutions to the ordinary 3-space 
dimensional wave equation, we have described exact 
teleportation of a time-domain multipole signal 
recorded at $r_1$ to another radial location $r_2$, 
thereby recovering the signal which would eventually 
reach $r_2$. We have focused on three issues, presenting 
new results for each.
First, we have described the structure of the exact 
convolution kernels which define such teleportation, and 
the large-$\ell$ (polar spherical harmonic index) 
cancellation errors associated with their exact 
sum-of-poles representations (in finite-precision arithmetic). 
In particular, we have given a slightly more detailed 
derivation of Greengard, Hagstrom, and Jiang's result 
for the large-$\ell$ asymptotics of
the residues. Second, we have described 
an accurate procedure, based on the work in Ref.~\cite{BFL},
for sphere-to-sphere signal propagation for the ordinary 
wave equation. The key ingredients of the procedure
are (i) accurate evaluation of frequency domain kernels 
through the algorithm described in section~\ref{sec:evaluation}
and (ii) rational approximation (i.e.~kernel compression)
of these kernels through the algorithm described in 
section~\ref{sec:compress}.
Lacking precise theoretical understanding of 
compressed teleportation kernels, we have modeled (from
numerical data) the relationships between the parameters 
determining these approximations.
As suggested in Ref.~\cite{GHJ}, the kernel representations 
derived in Ref.~\cite{GHJ} might be compressed using only 
the rational approximation step, thereby reducing the cost
of their evaluation.
The compressions reported here indicate that 
this strategy is possible, at least through $\ell = 256$.
Finally, we have focused on practical implementation,
demonstrating through simulations that 
teleportation can be performed accurately.

While this article has presented results only for the 
ordinary wave equation \eqref{eq:3d}, our methods can 
be extended to certain other hyperbolic PDEs, including 
those posed on non-flat geometries for which backscatter 
effects make the separation of "ingoing" and "outgoing" 
waves particularly vexing. Indeed, the main motivation 
for this work has been to compute asymptotic gravitational 
wave signals \cite{BFL}. Here the governing equations 
are those of Regge-Wheeler and Zerilli, respectively
describing axial and polar perturbations of spherically
symmetric blackholes. We are unaware of closed-form
time-domain representations for the asymptotically 
outgoing solutions to these equations. Moreover,
Laplace transforms of the corresponding (boundary and 
teleportation) kernels feature branch cuts in their
domains of analyticity \cite{Lau3}. Nevertheless, our approach
carries over to this more complicated gravitational
scenario. The idea is to again base kernel compression
on evaluation of the exact kernel through 
\eqref{eq:flatspacePhi} (other numerical
techniques are also required \cite{Lau1}). Based on 
this work, we conjecture that the 
methods described here can be applied to linear 
hyperbolic systems which (apart from inhomogeneities) are 
time independent and rotationally invariant. 
Extension of these methods to systems with
spheroidal invariance is an open problem.

Finally, we note that compressed kernels are available 
at \cite{KernelsWeb}. 
To date, we have mostly posted radiation boundary and
teleportation kernels for the Regge-Wheeler and Zerilli
equations. However, we intend to add the kernels
used in the numerical experiments documented here.

\section{Acknowledgements}
SRL gratefully acknowledges support from NSF grant 
No.~PHY 0855678 to the University of New Mexico, with
which infrastructure for our approximations was 
developed. SEF acknowledges support from the Joint Space Science 
Institute and NSF Grants No.~PHY 1208861 and No.~PHY 1005632 
to the University of Maryland, NSF Grants PHY-1306125 and 
AST-1333129 to Cornell University, and by a grant from 
the Sherman Fairchild Foundation. For insights and helpful 
comments we wish to thank Thomas Hagstrom and Akil Narayan.

\appendix 
\section{Residues of the frequency domain kernel} \label{sec:appendixA}
This first appendix derives the alternative
expression \eqref{eq:afromscaled}
for the residues $a_{\ell j}(r_1,r_2)$ in \eqref{eq:residues}.
We start with
\begin{align}
a_{\ell j}(r_1,r_2) & = 
\frac{W_\ell(b_{\ell j}r_2/r_1)}{r_1W_\ell'(b_{\ell j})} ,
\end{align}
a formula which agrees with \eqref{eq:residues}. We now
write $a_{\ell j}(r_1,r_2)$ in terms of standard 
special functions, starting with MacDonald's function.
From \eqref{eq:KfromW}, we get 
\begin{subequations}
\begin{align}
W_\ell(b_{\ell j} r_2/r_1) & = 
\sqrt{\frac{2 b_{\ell j}}{\pi}\frac{r_2}{r_1}}
e^{b_{\ell j} r_2/r_1} 
K_{\ell+1/2}(b_{\ell j} r_2/r_1)\\
W_\ell'(b_{\ell j}) & = 
\sqrt{\frac{2 b_{\ell j}}{\pi}}
e^{b_{\ell j}} K_{\ell+1/2}'(b_{\ell j}).
\end{align}
\end{subequations}
Note that $K_{\ell + 1/2}(z)$ is defined only on the slit plane
(due to the branch associated with the square root factor). 
Therefore, for odd $\ell$ the purely real root 
$b_{\ell,1+(\ell-1)/2}$ is not in the domain of analyticity.
Nevertheless, with appropriate cancellation of square-root
factors, the following expression is valid even for this root
\begin{equation}\label{eq:afromK} 
a_{\ell j}(r_1,r_2) = 
\frac{1}{r_1}\sqrt{\frac{r_2}{r_1}}e^{(r_2/r_1 - 1)b_{\ell j}}
\frac{K_{\ell+1/2}(b_{\ell j} r_2/r_1)}{K_{\ell+1/2}'(b_{\ell j})}.
\end{equation}

To use Olver's uniform asymptotic formulas for Bessel functions of 
large order and argument (see in particular AS {\bf 9.3.37} and
{\bf 9.3.45}), we need to express \eqref{eq:afromK} in terms of
the first Hankel function.
To this end, we start with the first equation in AS {\bf 10.2.15} 
(which involves the first spherical Hankel function defined in 
AS {\bf 10.1.16}):
\begin{align}\label{eq:K_h}
\sqrt{\frac{\pi}{2z}} K_{\ell + 1/2}(z) =
\frac{\mathrm{i}\pi}{2}
e^{\mathrm{i}(\ell + 1)\pi/2}
h_{\ell}^{(1)}(z e^{\mathrm{i}\pi/2}),
\qquad
-\pi < \arg z \leq \frac{1}{2}\pi.
\end{align}
Next, using the relationship between spherical and cylindrical 
Hankel functions given in AS {\bf 10.1.1}, we get
\begin{align}\label{eq:K_H}
\sqrt{\frac{\pi}{2z}} K_{\ell + 1/2}(z) 
= \frac{\mathrm{i}\pi}{2} 
e^{\mathrm{i}(2\ell + 1)\pi/4} 
\sqrt{\frac{\pi}{2z}}
H_{\ell+1/2}^{(1)}(z e^{\mathrm{i}\pi/2}),
\qquad
-\pi < \arg z \leq \frac{1}{2}\pi.
\end{align}
There are two branch cuts associated with the right-hand expression. 
The first is the usual cut along the negative real axis (in the 
$z$-plane) associated with the square root. The second results from 
the branch associated with $H_{\ell+1/2}^{(1)}(\bullet)$; due to the 
rotated argument this cut is the positive imaginary axis (in the 
$z$-plane). Across both cuts the right-hand expression in 
\eqref{eq:K_H} jumps by a sign, whereas, due to \eqref{eq:K_h} and 
the domain of analyticity for $h_{\ell}^{(1)}(\bullet)$, the 
left-hand expression in \eqref{eq:K_H} is analytic on the 
origin-punctured $z$-plane. Therefore, we work with the expression
\begin{equation}
K_{\ell + 1/2}(z) =
\epsilon
\frac{\mathrm{i}\pi}{2} e^{\mathrm{i}(2\ell + 1)\pi/4}
H_{\ell+1/2}^{(1)}(\mathrm{i}z),
\end{equation}
where $\epsilon = -1$ in the second quadrant and  $\epsilon = 1$
otherwise. Using this result, we cast \eqref{eq:afromK} into the 
form
\begin{equation}
a_{\ell j}(r_1,r_2) =
-\frac{\mathrm{i}}{r_1}\sqrt{\frac{r_2}{r_1}}
e^{(r_2/r_1 - 1)b_{\ell j}}
\frac{
H_{\ell+1/2}^{(1)}(\mathrm{i}b_{\ell j} r_2/r_1)
}{
H_{\ell+1/2}^{(1)\prime}(\mathrm{i}b_{\ell j})
}.
\end{equation}
The last expression leads directly to 
\eqref{eq:afromscaled} upon introduction of the
scaled zeros. In both the last expression and \eqref{eq:afromscaled}
the derivative in the denominator can be eliminated with the 
identity $2H_{\nu}^{(1)\prime}(z) =
H_{\nu-1}^{(1)}(z)-H_{\nu+1}^{(1)}(z)$.

\section{Conditioning and error bounds of teleportation 
with respect to data perturbation} \label{app}
We now consider conditioning of teleportation with respect 
to data perturbation and related error bounds.
In this appendix $(\ell,m)$ indices and, often, radial 
arguments $r_1,r_2$ on $\Psi$ and $\Phi$ are suppressed.
Let 
\begin{align}\label{eq:Aconv}
\Phi * \delta \Psi (t) = 
\int_0^t \Phi(t-t') 
\delta\Psi(t') dt', \qquad t \in [0,T].
\end{align}
Here $\delta\Psi(t',r_1) = \Psi(t',r_1) - \Psi_h(t',r_1)$ is 
the time-series error due to numerical discretization, 
and $T$ is the final time. 
One version of Young's convolution inequality yields
\begin{align}
\| \Phi * \delta \Psi \|_{L_\infty (0,T) } \leq 
\|\Phi\|_{L_\infty (0,T) } \cdot 
\|\delta \Psi\|_{L_1(0,T)},
\end{align}
which with \eqref{eq:td_teleportation} immediately gives
(upon replacing the radial arguments)
\begin{align} \label{eq:errbound}
\big\|\Psi(\cdot+(r_2-r_1),r_2) - 
      \Psi_h(\cdot+(r_2-r_1),r_2)\big\|_{L_\infty(0,T)}
\leq \|\Phi(\cdot,r_1,r_2)\|_{L_\infty (0,T) } \cdot 
\|\delta \Psi(\cdot,r_1)\|_{L_1(0,T)}
+ \| \delta \Psi(\cdot,r_1)\|_{L_{\infty}(0,T)}.
\end{align}
Numerical evidence suggests~\cite{GHJ} that 
$\left| \Phi_\ell(t,r_1,r_2)\right|$ is 
maximal at $t=0$. Furthermore, 
\begin{align} \label{eq:kernLinf}
\Phi_\ell(0,r_1,r_2) = \frac{1}{r_1} 
\left[\frac{\ell \left( \ell+1 \right)}{2} 
\left(\frac{r_1}{r_2} -1 \right) \right],
\end{align}
which follows from Lemma 1 of Ref.~\cite{GHJ} and 
the scaling relation \eqref{eq:scaling}.

\bibliographystyle{elsarticle-num}

\end{document}